\documentclass[11pt]{article}
\usepackage{mathrsfs}
\usepackage{upgreek}
\usepackage{float}
\usepackage{amssymb,amsmath,amsthm,amsfonts}
\usepackage[dvips]{graphicx}
\usepackage{epsfig, rotate}
\usepackage{pdfpages}

\usepackage[margin=1in]{geometry}

\def\dref#1{(\ref{#1})}

\newcommand{\etal}{\textit{et al.~}}
\numberwithin{equation}{section}

\newtheorem{theorem}{Theorem}[section]

\newtheorem{lemma}{Lemma}[section]

\def\C{\mathbb{C}}

\begin{document}
\thispagestyle{empty}

\noindent {\bf \Large Asymptotic Distributions for Likelihood Ratio
Tests \\
for the Equality of Covariance Matrices}

\vspace{20pt}
 \noindent{\bf Wenchuan Guo$^a$, Yongcheng Qi$^b$}

\vspace{10pt}

\noindent $^a$Biometrics \& Data Sciences, Bristol Myers Squibb,
3551 Lawrenceville Princeton, Lawrence Township, NJ 08543, USA.

\vspace{10pt}

\noindent $^b$Department of Mathematics and Statistics, University
of Minnesota Duluth, 1117 University Drive, Duluth, MN 55812, USA.

\date{\today}

\vspace{20pt}




\noindent{\bf Abstract.} Consider $k$ independent random samples
from $p$-dimensional multivariate normal distributions. We are
interested in the limiting distribution of the log-likelihood ratio
test statistics for testing for the equality of $k$ covariance
matrices. It is well known from classical multivariate statistics
that the limit is a chi-square distribution when $k$ and $p$ are
fixed integers. Jiang and Yang~\cite{JY13} and Jiang  and
Qi~\cite{JQ15} have obtained the central limit theorem for the
log-likelihood ratio test statistics when the dimensionality $p$
goes to infinity with the sample sizes.  In this paper, we derive
the central limit theorem when either $p$ or $k$ goes to infinity.
We also propose adjusted test statistics which can be well
approximated by chi-squared distributions regardless of values for
$p$ and $k$. Furthermore, we present numerical simulation results to
evaluate the performance of our adjusted test statistics and the
log-likelihood ratio statistics based on classical chi-square
approximation and the normal approximation.

\vspace{0.3cm}


\newpage

\section{Introduction}\label{intro}

High-dimensional data are increasingly common  in many applications
of statistics, encountered particularly in biological,
meteorological, and financial studies. If we have $n$ observations
from a $p$-dimensional population, it is natural to set up an
$n\times p$ data matrix. The analysis of such matrices includes
hypothesis testing, interval estimation and modeling, and is
elaborated in many classical textbooks of multivariate analysis such
as Anderson~\cite{anderson58}, Muirhead~\cite{muirhead82},  and
Eaton~\cite{eaton83}. In these books, the properties of these tests,
estimations or models are based on  the assumption that the data
dimensionality $p$ is a fixed constant or is negligible compared
with the sample size $n$. However, due to the explosive development
and improvement of computing power and data collection methods, the
data dimension is often high. The traditional assumption of small
$p$ large $n$ is no longer valid because the dimensionality $p$ can
be proportionally large compared with the sample size $n$. Some
recent results in this situation have been found in Bai
\etal~\cite{BJYZ}, Ledoit and Wolf~\cite{LM02}, Jiang
\etal~\cite{JJY12}, Jiang and Yang~\cite{JY13}, Jiang and
Qi~\cite{JQ15},  Dette and D\"ornemann~\cite{DD2020}, and references
therein.

Given that $k$ independent samples are observed and each of them is
from a $p$-dimensional multivariate normal distribution, we are
interested in testing the equality of $k$ covariance matrices of
multivariate normal distributions under the assumption that either
$p\rightarrow \infty$, or $k\rightarrow \infty$, or both. The
results in this direction can be very useful when one considers
multivariate analysis of variance problems for high dimensional data
that are resulted from a large number of treatments.

We denote $N_p(\mu, \Sigma)$ as the multivariate normal distribution
with dimensionality $p$, where $\mu$ is the mean vector and $\Sigma$
is the covariance matrix. Assume $k \ge 2$ is an integer. Let
$\{x_{i1},\cdots, x_{in_i}$\}, $1\le i  \le k$ be $k$ independent
random samples from $N_p(\mu_i, \Sigma_i)$, where the mean vectors
$\mu_i$ and covariance matrices $\Sigma_i$ are unknown.  Consider
the following hypothesis test:
\begin{equation} \label{test} H_0:
\Sigma_1=\cdots=\Sigma_k ~\text{ vs }~ H_a: H_0 \text{ is not true}.
\end{equation}

Define
\[
\overline x_i={1 \over
n_i}\sum_{j=1}^{n_i}x_{ij},~~~A_i=\sum_{j=1}^{n_i}(x_{ij}-\overline
x_i)(x_{ij}-\overline x_i)', \; 1\le i\le k,
\]
and
\[
A=A_1+\cdots+A_k,~ \text{ and }~n=n_1+\cdots+n_k.
\]
The likelihood ratio test for \dref{test}, derived by
Wilks~\cite{wilks32}, is as follows:
\[
\Lambda_n=\frac{\prod_{i=1}^k(|A_i|)^{n_i/2}}{(|A|)^{n/2}}\cdot
\frac{n^{np/2}}{\prod_{i=1}^kn_i^{n_ip/2}}.
\]
The details can be found in Section 8.2.2 of
Muirhead~\cite{muirhead82}. The test statistic $\Lambda_n$ is
non-degenerate only if every $A_i$ is of full rank, which means that
$p < n_i$ is required for all $i=1,\dots, k$, i.e, $p<\min_{1\le
i\le k} n_i$. However, the drawback of this likelihood-ratio test is
that it is biased when the sample sizes $n_1, \dots, n_k$ are not
all equal, which means that the probability of rejecting $H_0$ in
\eqref{test} when $H_0$ is true can be larger than the probability
of rejecting $H_0$ with $H_0$ is false.
 This phenomenon  was first noticed by Brown~\cite{brown39} with $p=1$.
 In what follows,  we use instead the following modified likelihood ratio test statistic $\Lambda_n^*$ suggested by Bartlett~\cite{bartlett37}:
 \begin{equation}
\label{Bartlett} \Lambda_n^*=\frac{\prod_{i=1}^k(|A_i|)^{(n_i-1)/2}}
{(|A|)^{(n-k)/2}}\cdot
\frac{(n-k)^{(n-k)p/2}}{\prod_{i=1}^k(n_i-1)^{(n_i-1)p/2}}.
\end{equation}
$\Lambda_n^*$ is derived by substituting the sample size $n_i$ by the degree of freedom $n_i-1$ and sum of $n_i$ by $n-k$. It is proved by Sugiura and Nagao~\cite{SN68}
that this modified likelihood ratio test is unbiased when $k=2$. A generalization to general integer $k$ is made in Perlman~\cite{perlman80}.

Assume $p$ and $k$ are fixed integers. Under the null hypothesis in \dref{test}, Box~\cite{box49} proved that $-2\rho \log\Lambda_n^*$ has a chi-square limit
\begin{equation}
\label{chiapprox}
-2\rho \log\Lambda_n^*\stackrel{d}{\to} \chi_f^2
\end{equation}
as $\min_{1\le i\le k}n_i \rightarrow \infty$, where
\begin{equation}\label{fn}
f={1 \over 2}p(p+1)(k-1),
\end{equation}
and
\[
\rho=1-\frac{2p^2+3p-1}{6(p+1)(k-1)(n-k)}(\sum_{i=1}^k{n-k \over n_i-1}-1).
\]
Note that $\rho$ converges to $1$ as $\min_{1\le i\le k}
n_i\to\infty$. The classic theory for likelihood ratio test
statistics states that  $-2\log\Lambda_n^*\stackrel{d}{\to}
\chi_f^2$. The correction factor $\rho$ is introduced in
\eqref{chiapprox} for a better chi-square approximation rate.

The likelihood ratio test of size $\alpha$ based on the asymptotic distribution of $-2\rho \log\Lambda_n^*$ is to reject $H_0$ \dref{test} when the sample is falling on the region
\begin{equation}\label{chi-R}
R_1=\{-2\rho \log\Lambda_n^*>c_f(\alpha)\},
\end{equation}
 where $c_f(\alpha)$ is the $\alpha$-level upper quantile of $\chi_f^2$ distribution.
  The details can be found in Muirhead~\cite{muirhead82}, p.308.

In recent studies, the assumption on the fixed dimensionality $p$ is
relaxed to change with the sample size $n_i$ so that the ratio
$p/n_i$  converges to a constant between 0 and 1, i.e.,
\begin{equation}\label{guo1}
\lim_{n \rightarrow \infty}{p \over n_i}=y_i \in (0,1).
\end{equation}
For instance, using the Random Matrix Theory method, the correction to the traditional likelihood ratio test for testing the equality of two covariance matrices is given in Bai \etal \cite{BJYZ} as
\[
L={|A|^{n_1/2} \cdot |B|^{n_2/2} \over |n_1A/n+n_2B/n|^{n/2}},
\]
where the two sample covariance matrices are denoted as $A$ and $B$, and a central limit theorem is established for $-2\log L$ under assumption \dref{guo1}.
An extension to case $y_i=1$ is also considered in Jiang, Jiang and Yang~\cite{JJY12}.  Furthermore, for likelihood ratio test for the
equality of $k$ covariance matrices, Jiang and Yang~\cite{JY13} and Jiang and Qi \cite{JQ15} have obtained the following central limit theorem (CLT)
\begin{equation}\label{CLT2}
\frac{\log\Lambda_n^*-a_n}{b_n}\stackrel{d}{\to}
N(0,1),
\end{equation}
where $a_n$ and $b_n>0$ are asymptotic mean and standard deviation
as the functions of $n$ and $p$.  In their work, $k\ge 2$ is a fixed
integer, the dimensionality $p$ can be proportionally  large to the
sample size $n_i$ or much smaller than $n_i$, that is,
$p/n_i\rightarrow y_i \in [0,1]$ and $n_i-1 \ge p$ for $1 \le i\le
k$.

The central limit theorem \eqref{CLT2} can also be true when both
$k$ and $p$ change with the sample sizes. Guo~\cite{guo2014} and
Dette and D\"ornemann~\cite{DD2020} establish the asymptotic
normality under very restrictive conditions on both $k$ and $p$. See
Remark 3 in Section~\ref{main} for more details.

The likelihood ratio statistic cannot be applied to the test of the
homogeneity of several covariance matrices if the dimensionality $p$
is larger than some of the sample sizes since the corresponding
sample covariance matrices are singular. To overcome this problem,
Schott~\cite{Schott07}, Srivastava and Yanagihara~\cite{SY10}
propose two test statistics based on trace of sample covariance
matrices. Subsequent papers such as Li and Chen~\cite{LC12}, Cai,
Liu and Xia~\cite{CLX13}, Yang and Pan~\cite{YP17}, and Zheng {\it
et al.}~\cite{ZLGY20} have developed various test statistics which
can be also applied to non-normal populations.

Our objective in this paper is to investigate the limiting
distribution of the log-likelihood ratio test statistic
$-2\log\Lambda^*_n$ for the multivariate normal distributions under
the minimal conditions. We have pointed out that the likelihood
ratio test exits only if $1\le p<\min_{1\le i\le k}n_i$.  In this
paper, we assume that both $p$ and $k$ can change with sample sizes,
$1\le p<\min_{1\le i\le k}n_i$, and $\min_{1\le i\le k}n_i$ tends to
infinity.  We will show in Theorem~\ref{thm1} that
$-2\log\Lambda^*_n$, after properly normalized, converges in
distribution to the standard normal distribution if $\max(p,
k)\to\infty$. Without imposing this additional condition,  we
propose an adjusted log-likelihood ratio test statistic (ALRT) in
\eqref{zn} in Section~\ref{main} and prove in
Theorem~\ref{thm-chisq} that the ALRT statistic is approximately
distributed a chi-square with $\frac12p(p+1)(k-1)$ degrees of
freedom. Neither  $-2\log\Lambda^*_n$ nor $-2\rho\log\Lambda^*_n$
shares this property  when $p$ or $k$ are too large.


The rest of the paper is organized as follows.  In
Section~\ref{main}, we give the main results in the paper, including
the CLT for $-2\log\Lambda^*_n$  and chi-square approximation to
the ALRT statistics and some discussions;
 In Section~\ref{sim}, we carry out a simulation study to compare the performance of $-2\log\Lambda^*_n$ based on the
 classical chi-square approximation, the normal approximation and the ALRT statistics based on the chi-square
 approximation. All proofs are provided in Section~\ref{proofs}.

\section{Main results}\label{main}

We first introduce some notations.  Let $\Gamma(x)$ denote the gamma
function given by
\[
\Gamma(x)=\int^\infty_0t^{x-1}e^{-t}dt,~~~x>0.
\]
Define the digamma function $\psi$ by
\begin{equation}\label{digammaf}
\psi(x)=\frac{d\log \Gamma(x)}{dx}=\frac{\Gamma'(x)}{\Gamma(x)}, ~~~x>0.
\end{equation}


Define
\begin{equation}\label{xi1}
\xi(x)=-2(\log(1-x)+x),~~~~~ x\in[0,1).
\end{equation}
 $\xi(x)$ is nonnegative in $[0,1)$.




Throughout this paper, we write $n=\sum^n_{i=1}n_i$ and assume that
the minimal sample size $\min_{1\le i\le k}n_i$ goes to infinity.
Furthermore, both the dimensionality parameter $p$ and number of the
populations $k$ can diverge. For this reason, we write $p=p_n$ and
$k=k_n$ to indicate that both parameters may change with the sample
sizes.

Define
\begin{eqnarray}\label{mun}
\mu_n&=&(n-k_n)\sum^{p_n}_{j=1}\psi(\frac{n-k_n+1-j}{2})-\sum_{i=1}^{k_n}(n_i-1)\sum^{p_n}_{j=1}\psi(\frac{n_i-j}{2})\\
&&~~+p_n\sum^{k_n}_{i=1}(n_i-1)\log(n_i-1)-p_n(n-k_n)\log(n-k_n)\nonumber
\end{eqnarray}
and
\begin{equation}\label{sigman}
\sigma_n^2=\sum^{k_n}_{i=1}(n_i-1)^2\xi(\frac{p_n}{n_i})-(n-k_n)^2\xi(\frac{p_n}{n-k_n+1})+p_n(k_n-1).
\end{equation}

\begin{theorem}\label{thm1} Assume $1\le p_n<\min_{1\le i\le k_n}n_i$ and $\min_{1\le i\le
k_n}n_i\to\infty$. Let $\Lambda_n^*$ be as in \dref{Bartlett}. Then,
under $H_0$ in \dref{test},
\begin{equation}\label{CLT}
\frac{-2\log \Lambda_n^*-\mu_n}{\sigma_n}\overset{d}\to N(0,1)~~~
\mbox{ as }n\to\infty
\end{equation}
if and only if $\max(p_n, k_n)\to\infty$ as $n\to\infty$.
\end{theorem}


\vspace{10pt}

Based on the Theorem~\ref{thm1}, the rejection region of the test
\dref{test} at level $\alpha$  is
\begin{equation}\label{clt-R}
R_2=\{(-2\log \Lambda_n^*-\mu_n)/\sigma_n >z_\alpha\},
\end{equation}
where $z_\alpha$ is $\alpha$ level critical value of the standard normal distribution.


The chi-square appropriation in \eqref{chiapprox} and the central
limit theorem in Theorem~\ref{thm1} exhibit two different types of
limiting distributions for inference. As our simulation study in
Section~\ref{sim} indicates, when $\max(k_n,p_n)$ is small, the
chi-square approximation \eqref{chiapprox} works very well; while
$\max(p_n,k_n)$ is large, the normal approximation prevails.  In
practice, one needs to decide which approximation should be used.
The two theoretical results don't clearly suggest when to use the
normal approximation and when to use the chi-square approximation.
It is desirable if we are able to construct a test statistic which
is a linear function of $-2\log\Lambda^*_n$ with a unified
asymptotic distribution regardless of values of $p_n$ and $k_n$. The
chi-square distributions are better choices than the normal
distribution for this purpose  since  $-2\log\Lambda^*_n$ is not
asymptotically normal when both $p$ and $k$ are small. We note that
the ratio of the mean and the variance of a chi-square distribution
is $1:2$. However, neither $-2\log\Lambda^*_n$ nor
$-2\rho\log\Lambda^*_n$ is well approximated by a chi-square
distribution when $\max(p_n, k_n)$ is too large as the ratio of
their asymptotic means and variances is not very close to $1:2$ in
this case.

We propose the following adjusted log-likelihood ratio test (ALRT) statistic
\begin{equation}\label{zn}
Z_n=(-2\log\Lambda_n^*)\sqrt{\frac{2f_n}{\sigma_n^2}}+f_n-\mu_n\sqrt{\frac{2f_n}{\sigma_n^2}},
\end{equation}
where $f_n=\frac12p_n(p_n+1)(k_n-1)$ as given in \eqref{fn}, $\mu_n$
and $\sigma_n^2$ are defined in \eqref{mun} and \eqref{sigman},
respectively. Note that the ALRT is essentially the LRT statistic
since it is a linear function of the log-likelihood ratio test
statistic, $-2\log\Lambda_n^*$.  This new test statistic achieves an
asymptotic mean $\frac{1}{2}p_n(p_n+1)(k_n-1)$ and an asymptotic
variance $p_n(p_n+1)(k_n-1)$. We have the following  theorem
regarding the chi-square approximation of $Z_n$ statistic.

\begin{theorem}\label{thm-chisq}
Assume $p_n\ge 1$ and $k_n\ge 2$ are two sequences of integers,
$p_n<\min_{1\le i\le k_n}n_i$ and $\min_{1\le i\le
k_n}n_i\to\infty$. Then we have under $H_0$ in \dref{test} that
\begin{equation}\label{chisquare}
\lim_{n\to\infty}\sup_{-\infty<x<\infty}|P(Z_n\le x)-P(\chi^2_{f_n}\le x)|=0.
\end{equation}
\end{theorem}

\vspace{20pt}

Based on Theorem~\ref{thm-chisq}, the rejection region of the test
\dref{test} at level $\alpha$ is
\begin{equation}\label{alrt-R}
R_3=\{Z_n>c_{f_n}(\alpha)\},
\end{equation}
where $c_{f_n}(\alpha)$ is the $\alpha$-level critical value of a chi-square distribution with $f_n=\frac12p_n(p_n+1)(k_n-1)$ degrees of freedom.

Under some restrictive conditions,  we can avoid using the digamma
function when we approach the mean of $-2\log\Lambda^*_n$.  We have
the following theorem.

\begin{theorem}\label{thm-special}
 Assume $p_n$ and $k_n\ge 2$ are two sequences of positive integers such that $\min_{1\le i\le k_n}n_i\to\infty$,
  $\sqrt{k_n}/\min_{1\le i\le k_n}n_i\to 0$ as $n\to\infty$, and there exists a constant $\delta\in (0,1)$
  such that $p_n\le \delta\min_{1\le i\le k_n}n_i$ for all large $n$.  Define
\begin{equation}\label{barmun}
\bar\mu_n=(n-k_n)(p_n-n+k_n+\frac12)\log(1-\frac{p_n}{n-k_n+1})-\sum^{k_n}_{i=1}(n_i-1)(p_n-n_i+\frac32)\log(1-\frac{p_n}{n_i})-p_n(k_n-1).
\end{equation}
\noindent (a). If $Z_n$ given in \eqref{zn} is redefined with
$\mu_n$ being replaced by $\bar\mu_n$ , then \eqref{chisquare} holds
under $H_0$ in
\eqref{test};\\
\noindent (b). Additionally, if $\max(p_n, k_n)\to\infty$ as
$n\to\infty$, then $(-2\log \Lambda_n^*-\bar\mu_n)/\sigma_n$
converges to $N(0,1)$ under $H_0$ in \eqref{test}.
\end{theorem}

\vspace{10pt}

\noindent \textbf{Remark 1.} In both Theorems~\ref{thm1} and
\ref{thm-chisq}, we assume $1\le p_n<\min_{1\le i\le k_n}n_i$ and
$\min_{1\le i\le k_n}n_i\to\infty$ as $n\to\infty$. The first
condition ensures the likelihood ratio test statistics exist and the
second condition guarantees that the minimal sample size goes to
infinity so that the asymptotic properties of $-2\log\Lambda^*_n$
can be investigated. Both conditions are necessary and cannot be
weakened. Theorem~\ref{thm1} shows that a necessary and sufficient
condition for the central limit theorems is
$\lim_{n\to\infty}\max(p_n,k_n)=\infty$, which is equivalent to
$\lim_{n\to\infty}f_n=\infty$, where
$f_n=\frac{1}{2}p_n(p_n+1)(k_n+1)$. If $f_n$ is bounded or fixed,
the distribution of $-2\log \Lambda_n^*$ is well approximated by
$\chi_{f_n}^2$ distribution. More generally, Theorem~\ref{thm-chisq}
ensures that the ALRT statistic, $Z_n$ as defined in \eqref{zn}, is
well approximately by $\chi_{f_n}^2$.

\vspace{10pt}

\noindent \textbf{Remark 2.}  In both
Theorems~\ref{thm1} and \ref{thm-chisq}, the digamma function $\psi$
as defined in \eqref{digammaf} is involved in the computation of the
asymptotic mean $\mu_n$ of $-2\log\Lambda_n^*$ via equation
\eqref{mun}. We notice that the digamma function $\psi$ is built in
software \textbf{R} and can be calculated using the \verb|digamma|
function in \textbf{R}.

\vspace{10pt}

\noindent \textbf{Remark 3.}  Assume the null hypothesis in
\eqref{test} holds, $1\le p_n<\min_{1\le i\le k_n}n_i$, and
$\min_{1\le i\le k_n}n_i\to\infty$.  Our Theorem~\ref{thm1}
indicates that $\lim_{n\to\infty}\max(p_n, k_n)=\infty$ is the
necessary and sufficient for the central limit limit theorems of
$-2\log \Lambda_n^*$. In Guo~\cite{guo2014}, it is shown that
\[
\frac{-2\log \Lambda_n^*-\bar\mu_n}{\sigma_n}\overset{d}\to
N(0,1)~~~ \mbox{ as }n\to\infty
\]
under additional conditions that $k_n/\min_{1\le i\le k_n}n_i\to 0$,
$\min(p_n, k_n)\to\infty$ as $n\to\infty$, and there exists a
constant $\delta\in (0,1)$ such that $p_n\le \delta\min_{1\le i\le
k_n}n_i$ for all large $n$. Theorem 3 in Dette and
D\"ornemann~\cite{DD2020} states that $-2\log \Lambda_n^*$, after
properly normalized, converges in distribution to the standard
normal under a different set of additional assumptions that
$p_n/(n-k_n)\to 0$, $\min(p_n, k_n)\to\infty$,
$k_n/\sqrt{p_n(n-k_n)}\to 0$ as $n\to\infty$, and there exists a
constant $\delta\in (0,1)$ such that $p_n\le \delta\min_{1\le i\le
k_n}n_i$ for all large $n$.   Clearly, either in Guo~\cite{guo2014}
or Dette and D\"ornemann~\cite{DD2020},  the following situations
are not under consideration:

\noindent (a) $\lim_{n\to\infty}p_n/\min_{1\le i\le k_n}n_i=1$;

\noindent (b) $\lim_{n\to\infty}k_n/(p_n\max_{1\le i\le k_n}n_i)>0$.

\vspace{10pt}

\noindent \textbf{Remark 4.} As we have mentioned earlier,
conditions $1\le p_n<\min_{1\le i\le k_n}n_i$ and $\min_{1\le i\le
k_n}n_i\to\infty$ are necessary for discussion on the asymptotic
distribution of $-2\log\Lambda^*_n$. Theorem~\ref{thm-chisq}
indicates that the new test statistic $Z_n$ can be approximated by
chi-square distribution uniformly over $1\le p_n<\min_{1\le i\le
k_n}n_i$ and $k_n\ge 2$ as $\min_{1\le i\le k_n}n_i\to\infty$.
However,  there is a transition from chi-square to the normal in the
limiting distribution of $-2\log\Lambda^*_n$.




\section{Simulation study}\label{sim}

In this section, we  carry out some simulation studies to compare
the performance for the likelihood ratio test statistics based on
three different approximations, including the classic chi-square
approach \eqref{chiapprox}, the CLT approach \eqref{CLT} and ALRT
approach \eqref{chisquare}.  The corresponding rejection regions for
a size $\alpha$ test on \eqref{test} are given in equations
\eqref{chi-R}, \eqref{clt-R} and \eqref{alrt-R}, respectively.

We compare the three test statistics as follows:  a).  We plot the
histograms of the three statistics and see how well the histograms
for each statistic fits its theoretical distribution under the null
hypothesis of \eqref{test}; b).  We estimate the sizes (or type I
errors) and determine which test statistic is more accurate; c). We
estimate the power for each test statistic under some alternative
hypotheses of \eqref{test}.

For each combination of $p$ and $k$ with $p=3, 20, 50, 95$ and $k=5,
20, 50$, we set $n_i=100$ for $1 \le i \le k$.  We generate $k$
random samples with sample sizes $n_1, \cdots, n_k$ respectively
from multivariate normal distributions with selected covariance
matrices, and then we calculate three test statistics. Histograms,
estimation of sizes and powers are based on $10,000$ replicates. In
our histograms and tables for sizes and powers for three test
statistics,  we use ``Chisq", ``CLT" and ``ALRT" to denote the
classic chi-square approximation \eqref{chiapprox}, normal
approximation \eqref{CLT}, and adjusted chi-square approximation
\eqref{chisquare}, respectively.

\subsection{Histograms and empirical sizes}

We generate random samples from $k$ multivariate normal
distributions under the null hypothesis of \eqref{test}. Denote the
common variance matrix as $\Sigma$ which is a $p$ by $p$ matrix.
Since the distribution function of $\Lambda^*_n$ doesn't depend on
the matrix $\Sigma$, we can simply set $\Sigma=I_p$ as a $p$ by $p$
identity matrix.

For each combination of $p$ and $k$, we plot the histograms for the
three standardized test statistics, namely $-2\rho\log\Lambda_n^*$,
$(-2\log\Lambda_n^*-\mu_n)/\sigma_n$, and $Z_n$. The histograms for
these statistics are also compared with their approximate density
functions -- a chi-squared distribution with $\frac12p(p+1(k-1)$
degrees of freedom, the standard normal distribution, and a
chi-squared distribution with $f=\frac12p(p+1)(k-1)$ degrees of
freedom, respectively, as established in \eqref{chiapprox},
\eqref{CLT} and \eqref{chisquare}. See Figure~\ref{histogram-k3} and
Figures S2 and S3 in the supplement.

The histograms in Figure~\ref{histogram-k3} and Figures S2 and S3
reveal the overall performance for the three test statistics. We can
conclude that the classic chi-square approximation works well when
$p$ is small and gets worse as $p$ increases, the normal
approximation starts to get better as $f=\frac{1}{2}p(p+1)(k-1)$
increases, and only the histograms for the adjusted test statistic
fit the chi-squared curves in all cases.

The empirical sizes for all three tests at level $\alpha=0.05$ are
reported in Table~\ref{table1}. A good test in terms of the type I
errors (or sizes) will have an actual size close to its nominal
level $0.05$. The adjusted test statistic outperforms over the other
two test statistics since its sizes are close to $0.05$ in all
cases. In particular, the classic chi-square approximation performs
well only when $p$ is relatively small. When $p$ is large, the null
hypothesis is rejected with a probability close to one.

\begin{figure}[H]
\centering
\includegraphics[height=2.0in,width=6in]{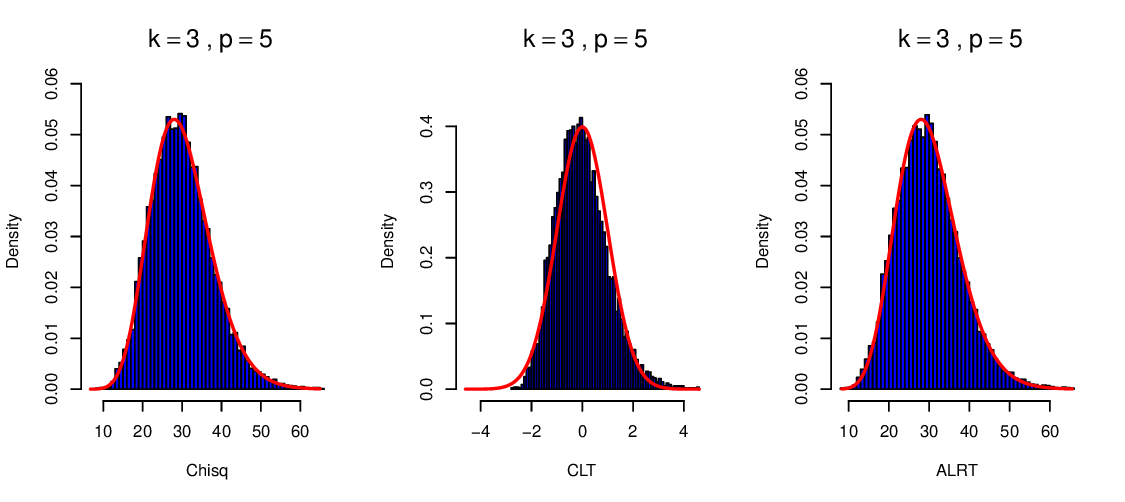}
\includegraphics[height=2.0in,width=6in]{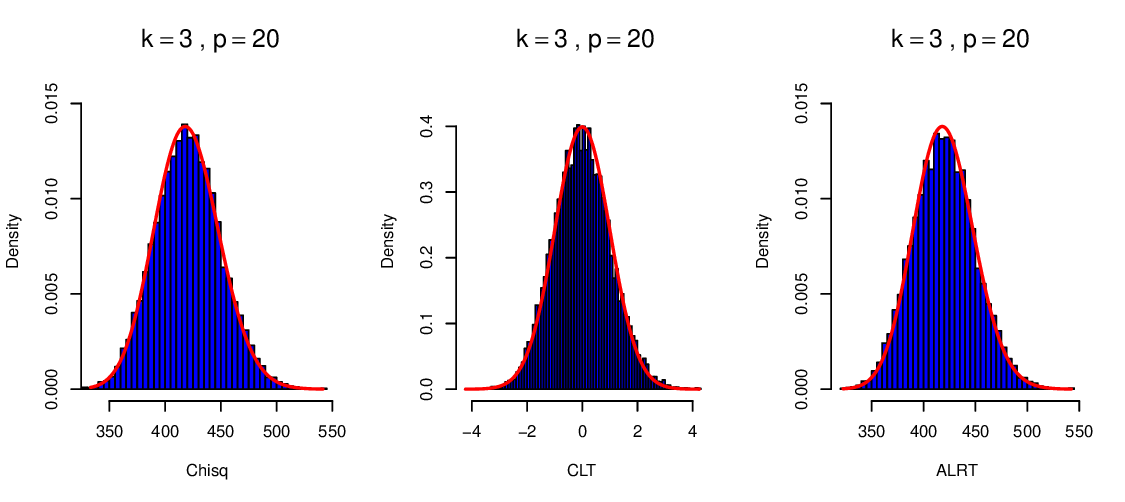}
\includegraphics[height=2.0in,width=6in]{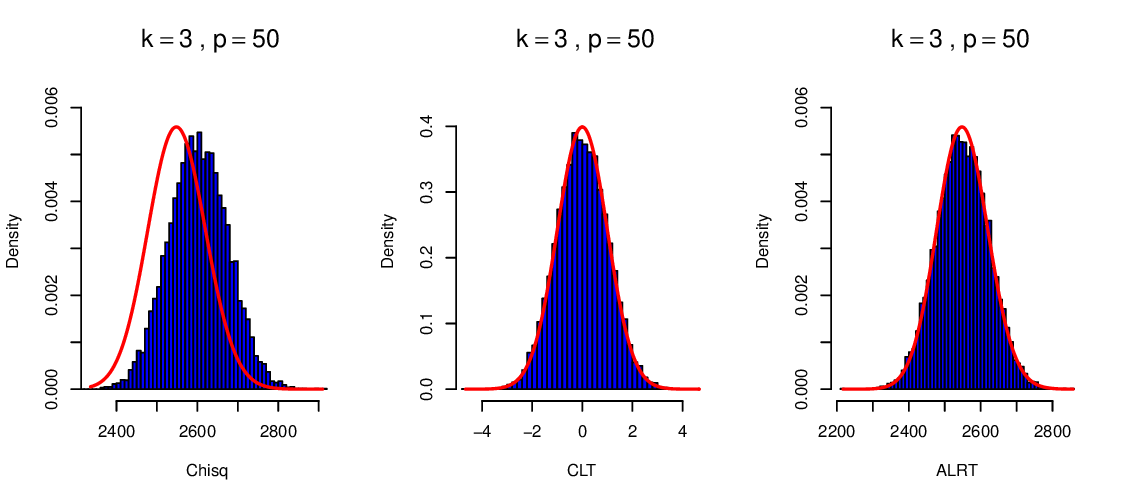}
\includegraphics[height=2.0in,width=6in]{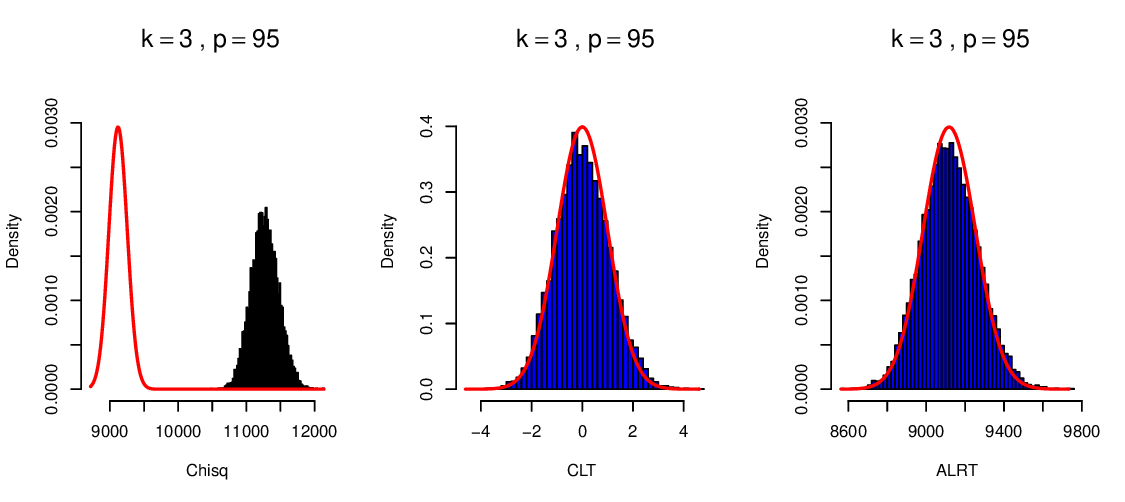}
\caption{Histograms of $-2\rho \log\Lambda_n^*$ (Chisq),
$(-2\log\Lambda_n^*-\mu_n)/\sigma_n$ (CLT), and  $Z_n$ (ALRT), where
$n_i=100$, $1 \le i\le k$, with $k=3$ and $p=5,20,50,95$.}
 \label{histogram-k3}
\end{figure}

\begin{table}[ht]
\caption{The size of LRT for equality of $k$ covariance
matrices}\label{table1}
 \centering
\begin{tabular}{l c c c}
\hline
Conditions &Chisq &CLT & ALRT  \\
\hline
 $n_i=100, ~k=~3, ~p=~5$ & 0.0501 &0.0636 &0.0529 \\
 $n_i=100, ~k=~3, ~p=20$ & 0.0575 &0.0608 &0.0569 \\
 $n_i=100, ~k=~3, ~p=50$ & 0.1861 &0.0539 &0.0521 \\
 $n_i=100, ~k=~3, ~p=95$ & 1.0000 &0.0689 &0.0678 \\
 $n_i=100, ~k=20, ~p=~5$ & 0.0511 &0.0576 &0.0529 \\
 $n_i=100, ~k=20, ~p=20$ & 0.0639 &0.0542 &0.0530 \\
 $n_i=100, ~k=20, ~p=50$ & 0.8135 &0.0533 &0.0529 \\
 $n_i=100, ~k=20, ~p=95$ & 1.0000 &0.0618 &0.0614 \\
 $n_i=100, ~k=50, ~p=~5$ & 0.0490 &0.0540 &0.0505 \\
 $n_i=100, ~k=50, ~p=20$ & 0.0755 &0.0525 &0.0519 \\
 $n_i=100, ~k=50, ~p=50$ & 0.9937 &0.0554 &0.0551 \\
 $n_i=100, ~k=50, ~p=95$ & 1.0000 &0.0619 &0.0618 \\
 \hline
\end{tabular}
\end{table}

\subsection{Comparison of powers}\label{esti_power}

The power of a test is the probability of rejecting $H_0$ when $H_0$
is false. For each combination of $k$ and $p$, we take $k$ random
samples, and the $i$-th sample is  from a multivariate normal
distribution with mean vector $\textbf{0}$ and covariance matrix
$\Sigma_i=I_p+\frac{0.05(i-1)}{\sqrt{k}}J_p$ for $1\le i\le k$,
where $I_p$ is a $p$ by $p$ identity matrix, and $J_p$ is a $p$ by
$p$ matrix with all entries being equal to $1$. Again, all sample
sizes are the same with $n_1=\cdots=n_k=100$ in the study. The
powers for the three test statistics are estimated based on $10,000$
replicates, and they are reported in Table~\ref{table2}.

We first note that all three test statistics are linear combinations
of $-2\log\Lambda_n^*$, the corresponding rejection regions at level
$\alpha$ are of the same type as $\{-2\log\Lambda_n^*>c\}$. The only
difference is that the values of cutoff $c$ corresponding to the
three approximations are different.  That means the three rejections
$R_1, R_2$ and $R_3$, defined in \eqref{chi-R}, \eqref{clt-R} and
\eqref{alrt-R}, respectively,  are nested each other.  A large
rejection region has a larger size of test and a larger power as
well.  From Table~\ref{table2}, we can see that whenever any two
test statistics have really close sizes, they have similar powers.
When the classic chi-square approximation fails, that is, it
produces a large size, it also yields a large power.  However, its
large power does not make any sense here.  So we compare the normal
approximation and the chi-square approximation based on the test
statistic $Z_n$.  We observe that both test statistics have
comparable powers when the normal approximation achieves a
reasonable size.

\begin{figure}[H]
\centering
\includegraphics[height=3.5in,width=3in]{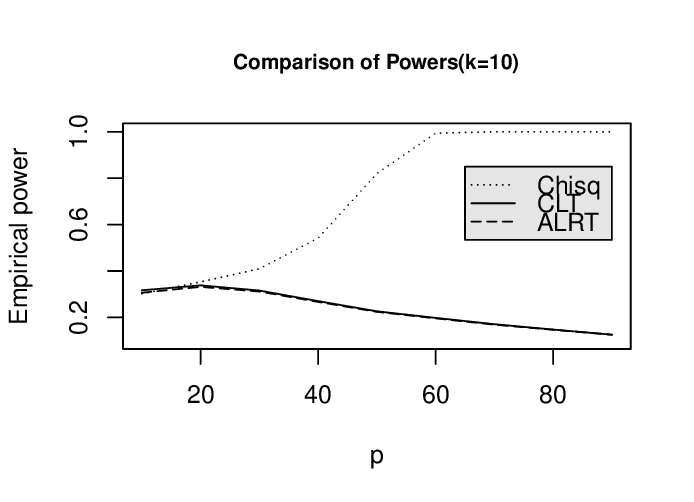}
\includegraphics[height=3.5in,width=3in]{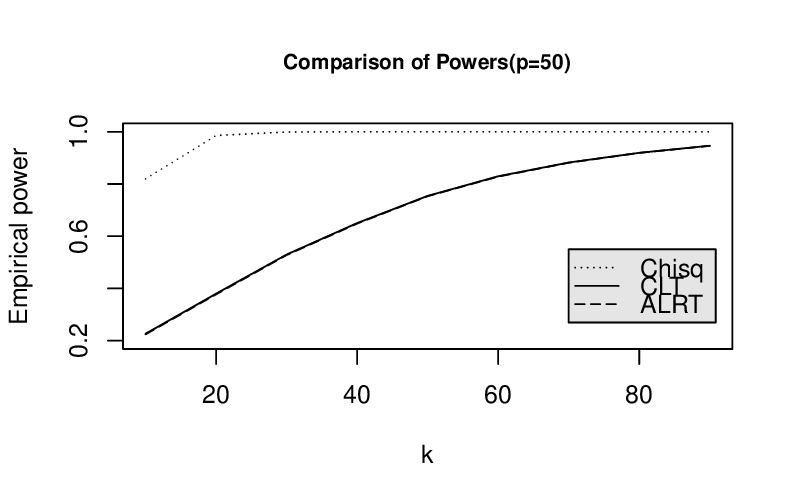}
\caption{Plots of powers for $-2\rho \log\Lambda_n^*$ (Chisq),
$(-2\log\Lambda_n^*-\mu_n)/\sigma_n$ (CLT), and  $Z_n$ (ALRT)
 with $n_1=\cdots=n_{k}=100$. The left one shows how the powers change with $p$ when $k=10$,
 and the right one displays how the powers change with $k$ when $p=50$.}
 \label{powercurve}
\end{figure}

Finally, we will demonstrate how the powers of the three test
statistics change with $p$ or $k$ while other variables remain
constant. For illustration purpose, we set $k=10$ and $n_i=100$ when
$p$ increases from $10$ to $90$, and set $p=50$ and $n_i=100$ when
$k$ increases from $10$ to $90$.  The covariance structures are the
same as specified in the beginning of this subsection. The curves of
the powers for the three test statistics are plotted in
Figure~\ref{powercurve}. We see that the power of the classic
chi-square statistic (Chisq) is much larger when $p$ is bigger, and
this is the case for the two plots in the figure. Keep in mind that
we have set $p=50$ in the second plot in Figure~\ref{powercurve}.
This phenomenon happens since the classic chi-square approach
results in a much larger type I error than the nominal level $0.05$
as in Table~\ref{table1}.  Two other methods are very close in terms
of power.

\begin{table}[ht]
\caption{The power of LRT for equality of $k$ covariance
matrices}\label{table2}
 \centering
\begin{tabular}{l c c c}
\hline
Conditions &Chisq &CLT & ALRT  \\
\hline
$n_i=100, ~k=~3, ~p=~5$ & 0.0744 &0.0958 &0.0778 \\
$n_i=100, ~k=~3, ~p=20$ & 0.1203 &0.1260 &0.1194 \\
$n_i=100, ~k=~3, ~p=50$ & 0.3249 &0.1233 &0.1197 \\
$n_i=100, ~k=~3, ~p=95$ & 1.0000 &0.0987 &0.0976 \\
$n_i=100, ~k=20, ~p=~5$ & 0.5587 &0.5809 &0.5656 \\
$n_i=100, ~k=20, ~p=20$ & 0.7025 &0.6728 &0.6690 \\
$n_i=100, ~k=20, ~p=50$ & 0.9858 &0.3798 &0.3775 \\
$n_i=100, ~k=20, ~p=95$ & 1.0000 &0.1478 &0.1473 \\
$n_i=100, ~k=50, ~p=~5$ & 0.9992 &0.9992 &0.9992 \\
$n_i=100, ~k=50, ~p=20$ & 0.9971 &0.9951 &0.9951 \\
$n_i=100, ~k=50, ~p=50$ & 1.0000 &0.7544 &0.7534 \\
$n_i=100, ~k=50, ~p=95$ & 1.0000 &0.2284 &0.2278 \\
\hline
\end{tabular}
\end{table}

\subsection{Concluding remark}

In this paper, we study the limiting distribution of the likelihood
ratio test statistics for testing equality of the covariance
matrices from $k$ multivariate normal populations. We allow both $k$
and the dimensionality $p$ of the populations to change with the
sample sizes. We establish the central limit theorem for the
log-likelihood test when $\max(k, p)$ tends to infinity. Further, we
propose an adjusted log likelihood test statistic and prove that the
distribution of the new test statistic can be well approximated by a
chi-squared distribution regardless of how $k$ and $p$ change with
the sample sizes.  This new test outperforms in terms of accuracy in
type I error and has the comparable power as the classic chi-square
approximation or the normal approximation when they are applicable.

\section{Proofs}\label{proofs}

The standard notation to  describe the limiting behavior  of a
function is denoted by big $O$ and small $o$. Precisely, for two
sequences of numbers $\{a_n, n\ge 1\}$ and $\{b_n, n\ge1\}$,
$a_n=O(b_n)$ as $n \rightarrow \infty$ means that
$\lim\sup_{n\to\infty} |a_n/b_n| < \infty$ and $a_n=o(b_n)$ as $n
\rightarrow \infty$ means that $\lim_{n\to\infty}a_n /b_n=0$. The
symbol $a_n\sim b_n$ stands for $\lim_{n\to\infty}a_n/b_n=1$. For
two functions $f$ and $g$, we also use $f(x)=O(g(x)), f(x)=o(g(x))$,
and $f(x)\sim g(x)$ with the similar meaning as $x\rightarrow x_0
\in [-\infty,\infty]$.

The multivariate gamma function, denoted by $\Gamma_p(z)$, is
defined as
\begin{equation}\label{gammap}
\Gamma_q(z)=\pi^{q(q-1)/4}\prod_{i=1}^{q}\Gamma\Big(z-{1 \over 2}(i-1)\Big),
\end{equation}
for complex numbers $z$ with Re$(z)>{1 \over 2}(q-1)$, where $q\ge
1$ is an integer.  In our applications we only need to consider
$\Gamma_q(z)$ for real variable $z>\frac{q-1}{2}$. We will put much
effort in expand the function $\Gamma_q(z)$ below.

The function  $\xi(x)$ defined in \eqref{xi1} is nonnegative in $x\in [0,1 )$.
It is easy to see that $\xi(0)=0$ and $\xi'(x)=\frac{2x}{1-x}$. We have
\[
\xi(x)
=\xi(x)-\xi(0)\\
=\int_{0}^{x}\xi'(t)dt\\
=2\int_{0}^{x}\frac{t}{1-t}dt.
\]
By writing $t=ux$, we have $dt=xdu$, and
\[
\xi(x)=2\int_{0}^1 \frac{ux^2}{1-ux}du
=2x^2\int_{0}^{1}\frac{u}{1-ux}du,  ~~~~~x\in [0,1).
\]

Now we define
\begin{equation}\label{eta}
\eta(x)=\frac{\xi(x)}{x^2}=2\int_{0}^{1}\frac{u}{1-ux}du, ~~~~x\in [0,1).
\end{equation}
We can conclude that
\begin{equation}\label{eta-property}
\eta(0)=1, ~~\lim_{x\uparrow 1}\eta(x)=\infty, \mbox{ and }\eta(x)\ge 1 \mbox{ is increasing in } [0,1).
\end{equation}

\subsection{Some auxiliary lemmas}

\begin{lemma}\label{eta1} Assume $\{r_m, ~m\ge 2\}$ is any sequence of integers such that $1\le r_m<m$ and $r_m/m\to$ as $m\to\infty$. Then for any $\delta>0$
\[
\max_{r_m\le q<m}\big(\frac{m}{q(m-q)}\big)^{\delta}\frac{1}{\eta(\frac{q}{m})}\to 0~~~\mbox{as }m\to\infty,
\]
where $\eta$ is defined in \eqref{eta}.
\end{lemma}

\noindent{\it Proof.} One can easily verify
\begin{equation}\label{bound}
\frac{m}{q(m-q)}\le \frac{m}{m-1}\le 2, ~~~~~~~\mbox{ for } 1\le q<m
\end{equation}
and
\begin{equation}\label{bound1}
\min_{r_m\le q\le m-r_m}\frac{m}{q(m-q)}\le \frac{m}{r_m(m-r_m)}\to 0~~\mbox{ as }m\to\infty.
\end{equation}
Then it follows from equations \eqref{eta-property}, \eqref{bound} and \eqref{bound1} that
\begin{eqnarray*}
\max_{r_m\le
q<m}\big(\frac{m}{q(m-q)}\big)^{\delta}\frac{1}{\eta(\frac{q}{m})}
&=&\max\Big(\max_{r_m\le q\le
m-r_m}\big(\frac{m}{q(m-q)}\big)^{\delta}\frac{1}{\eta(\frac{q}{m})},
\max_{m-r_m\le q\le m}\big(\frac{m}{q(m-q)}\big)^{\delta}\frac{1}{\eta(\frac{q}{m})}\Big)\\
&\le& \max\Big(\big(\max_{r_m\le q\le m-r_m}\frac{m}{q(m-q)}\big)^\delta, \frac{2^{\delta}}{\eta(\frac{m-r_m}{m})}\Big)\\
&\to & 0
\end{eqnarray*}
as $m\to\infty$.  This completes the proof.  \hfill $\blacksquare$

\begin{lemma}\label{oxi} As $m\to\infty$
\begin{equation}\label{suff}
\max_{1\le q<m}\frac{q^2}{m^2(m-q)}\frac{1}{\xi(\frac{q}{m})}=o(1).
\end{equation}
and
\begin{equation}\label{suff2}
\max_{1\le q<m}\frac{q}{m(m-q)^2}\frac{1}{\xi(\frac{q}{m})}=o(1).
\end{equation}
\end{lemma}

\noindent{\it Proof.} Let $r_m$ be any sequence of integers satisfying Lemma~\ref{eta1}.
 Since
\[
\frac{q^2}{m^2(m-q)}\frac{1}{\xi(\frac{q}{m})}=\frac{1}{m-q}\frac{1}{\eta(\frac{q}{m})}
\]
and $\eta(x)\ge 1$, we have
\begin{eqnarray*}
\max_{1\le q<m}\frac{q^2}{m^2(m-q)}\frac{1}{\xi(\frac{q}{m})}&=&\max\Big(\max_{r_m\le q<m}\frac{q}{m}\frac{m}{q(m-q)}\frac{1}{\eta(\frac{q}{m})}, \max_{1\le q< r_m}\frac{1}{m-q}\frac{1}{\eta(\frac{q}{m})}\Big)\\
&\le &\max\Big(\max_{r_m\le q<m}\frac{m}{q(m-q)}\frac{1}{\eta(\frac{q}{m})}, \max_{1\le q<r_m}\frac{1}{m-q}\Big)\\
&=&\max\Big(\max_{r_m\le q<m}\frac{m}{q(m-q)}\frac{1}{\eta(\frac{q}{m})}, \frac{1}{m-r_m+1}\Big),
\end{eqnarray*}
which converges to zero as $m\to\infty$ according to Lemma~\ref{eta1}. This proves \eqref{suff}

\eqref{suff2} can be proved in the same manner and the details are omitted.  \hfill$\blacksquare$

\begin{lemma}\label{digamma} For any $x_0>0$, there exist constants $C_1$ and $C_2$ such that
\begin{equation}\label{gamma0}
|\psi(x)-\log x+\frac{1}{2x}+\frac{1}{12x^2}|\le \frac{C_0}{x^4},~~~x\ge x_0;
\end{equation}
\begin{equation}\label{gamma1}
|\psi'(x)-\frac1x-\frac{1}{2x^2}|\le \frac{C_1}{x^3},~~~x\ge x_0.
\end{equation}
\end{lemma}

\noindent{\it Proof.} Without loss of generality, assume $x_0\in (0,1)$.

It follows from Formulas 6.3.18 and 6.4.12
in Abramowitz and Stegun~\cite{Abramowitz1972} that
\[
\psi(x)=\log x-\frac{1}{2x}-\frac{1}{12x^2}+O(\frac{1}{x^4}),
\]
\[
\psi'(x)=\frac1x+\frac{1}{2x^2}+O(\frac{1}{x^3})
\]
as $x\to\infty$.  Therefore,  there exist constant $x_1\ge 1$ and $C>0$ such that \eqref{gamma0} and \eqref{gamma1}
hold for $x\ge C$ with $C_0=C_1=C$.  It remains to show that \eqref{gamma0} and \eqref{gamma1}
hold for $x_0\le x\le x_1$.  We show \eqref{gamma0} only for illustration. Since the function
$\psi(x)-\log x+\frac{1}{2x}+\frac{1}{12x^2}$ is continuous in $[x_0,x_1]$, we have $C_3:=\sup_{x\in [x_,x_1]}|\psi(x)-\log x+\frac{1}{2x}+\frac{1}{12x^2}|<\infty$.
By setting $C_0=\max(C, x_1^4C_3)$, we have
\[
|\psi(x)-\log x+\frac{1}{2x}+\frac{1}{12x^2}|\le C_3\le\frac{C_3x_1^4}{x^4}=\frac{C_0}{x^4}, ~~x\in [x_0,x_1].
\]
Obviously, we also have
\[
|\psi(x)-\log x+\frac{1}{2x}+\frac{1}{12x^2}|\le\frac{C}{x^4} \le \frac{C_0}{x^4}, ~~x\in [x_1,\infty).
\]
This proves \eqref{gamma0}.   The proof of Lemma~\ref{digamma} is done.  \hfill$\blacksquare$

\begin{lemma}\label{loggamma}
Define
\[
s(x)=(x-\frac12)\log(x)-x,~~x>0
\]
and
\begin{equation}\label{h(x)}
h(x)=\log \Gamma(x)-s(x), ~~x>0.
\end{equation}
Then we have
\begin{equation}\label{s-diff}
s'(x)=\log x-\frac{1}{2x}, ~ s''(x)=\frac{1}{x}+\frac{1}{2x^2}, ~
s'''(x)=-\frac{1}{x^2}-\frac{1}{x^3},
\end{equation}
and there exists a constant $c>0$ such that
\[
|h''(x)|\le \frac{C}{x^3},~~~ x\ge 1/4.
\]
\end{lemma}

\noindent{\it Proof.} \eqref{s-diff} is obvious. The last estimate in the lemma follows from \eqref{gamma1} by noting that $h''(x)=\psi'(x)-s''(x)$.
 \hfill$\blacksquare$

\begin{lemma}\label{x-delta} We have for any $\delta\ge 1$
\begin{equation}\label{myest}
0\le \frac{1}{x^\delta}-\frac{1}{y^\delta} \le \frac{\delta(y-x)}{x^{\delta}y},~~~y>x>0.
\end{equation}
\end{lemma}

\noindent{\it Proof.} In fact, for any $\delta\ge 1$, $y>x>0$
\[
0\le \frac{1}{x^\delta}-\frac{1}{y^\delta}=\delta\int^y_x\frac{1}{t^{1+\delta}}dt\le\frac{\delta}{x^{\delta-1}}\int^y_x\frac{1}{t^2}dt
=\frac{\delta}{x^{\delta-1}}(\frac{1}{x}-\frac{1}{y})= \frac{\delta(y-x)}{x^{\delta}y},
\]
proving \eqref{myest}.
\hfill$\blacksquare$

\begin{lemma}\label{approximation12} As  $m\to\infty$,
\begin{equation}\label{var}
\sum^q_{i=1}(\frac1{m-i}-\frac{1}{m-1})=-\frac{q}{m}-\log(1-\frac{q}{m})+O(\frac{q}{m(m-q)}),
\end{equation}
\begin{equation}\label{mean}
\sum^q_{i=1}\Big(\log(m-i)-\log(m-1)\Big)=(q-m+\frac12)\log(1-\frac{q}{m})-\frac{(m-1)q}{m}+O(\frac{q}{m(m-q)}),
\end{equation}
\begin{equation}\label{m-i-m}
\sum_{i=1}^q({1 \over m-i}-\frac{1}{m-1})=-{q \over m}-\log(1-{q
\over m})-\frac{q}{2(m-1)^2}+o(\xi(\frac{q}{m})),
\end{equation}
\begin{equation}\label{sum2}
\sum^q_{i=1}\frac{1}{(m-i)^2}=\frac{q}{(m-1)^2}+o(\xi(\frac{q}{m}))
\end{equation}
and for any fixed $\delta\ge 2$
\begin{equation}\label{sum-delta}
\sum^q_{i=1}\frac{1}{(m-i)^{\delta}}= O(\frac{q}{(m-q)^{\delta-1}m})
\end{equation}
uniformly over $1\le q< m$ as $m\to\infty$.
\end{lemma}

\noindent {\it Proof.}  \eqref{var} and \eqref{mean} are proved in Lemma 5 in Qi, Wang and Zhang~\cite{QWZ19}.

Recall the relation of $\xi(x)$ and $\eta(x)$ from \eqref{eta} that $\xi(x)=x^2\eta(x)$ for $0<x<1$. Then we have
\[
\frac{q}{m(m-q)}\frac{1}{\xi(\frac{q}{m})}=\frac{m}{q(m-q)}\frac{1}{\eta(\frac{q}{m})},~~1\le q<m.
\]
Let $r_m$ be any sequence of integers such that $1\le r_m<m$ and $r_m/m\to 0$ as $m\to\infty$,  From from Lemma~\ref{eta1} we have
\[
\max_{r_m\le q<m}\frac{q}{m(m-q)}\frac{1}{\xi(\frac{q}{m})}\to 0~~~\mbox{as }m\to\infty
\]
and
\[
\max_{r_m\le q<m}\frac{q}{(m-1)^2}\frac{1}{\xi(\frac{q}{m})}\le \max_{r_m\le q<m}\frac{q}{m(m-q)}\frac{1}{\xi(\frac{q}{m})}\to 0 ~~~\mbox{as }m\to\infty,
\]
which combined with \eqref{var} yields that \eqref{m-i-m} holds uniformly over $r_m\le q<m$ as $m\to\infty$. It remains to show that \eqref{m-i-m} holds uniformly over $1\le q<r_m$ as $m\to\infty$.  In fact, by uisng Taylor's expansion we have  uniformly over  $1\le q<r_m$
\begin{eqnarray*}
\sum_{i=1}^q({1 \over m-i}-\frac{1}{m-1})&=&\sum_{i=1}^q\frac{1}{m-1}(\frac{1}{1-\frac{i-1}{m-1}}-1)\\
&=&\sum_{i=1}^q\frac{1}{m-1}\big(\frac{i-1}{m-1}+O((\frac{i-1}{m-1})^2)\big)\\
&=&\sum_{i=1}^q\big(\frac{i-1}{(m-1)^2}+O((\frac{(i-1)^2}{(m-1)^2})\big)\\
&=&\frac{q(q-1)}{2(m-1)^2}+O(\frac{q^3}{m^3})\\
&=&\frac{q^2}{2(m-1)^2}-\frac{q}{2(m-1)^2}+O(\frac{q^3}{m^3})\\
\end{eqnarray*}
and
\[
-\frac{q}{m}-\log(1-\frac{q}{m})=\frac{q^2}{2m^2}+O(\frac{q^3}{m^3}),
\]
which imply
\begin{eqnarray*}
\sum^q_{i=1}(\frac1{m-i}-\frac{1}{m-1})&=&-\frac{q}{m}-\log(1-\frac{q}{m})-\frac{q}{2(m-1)^2}+\big(\frac{q^2}{(m-1)^2}-\frac{q^2}{m^2}\big)+O(\frac{q^3}{m^3})\\
&=&-\frac{q}{m}-\log(1-\frac{q}{m})-\frac{q}{2(m-1)^2}+O(\frac{q^3}{m^3})\\
&=&-\frac{q}{m}-\log(1-\frac{q}{m})-\frac{q}{2(m-1)^2}+O(\frac{r_m}{m}\frac{q^2}{m^2})\\
&=&-\frac{q}{m}-\log(1-\frac{q}{m})-\frac{q}{2(m-1)^2r}+o(\xi(\frac{q}{m}))
\end{eqnarray*}
uniformly over $1\le q<r_m$ as $m\to\infty$. In the las step we use that $r_m/m\to 0$ as $m\to\infty$ and $\xi(\frac{q}{m})=\frac{q^2}{m^2}\eta(\frac{q}{m})\ge \frac{q^2}{m^2}$. This completes the proof of \eqref{m-i-m}.

For $1\le i\le q<m$, we have
\begin{eqnarray*}
0\le \frac{1}{(m-i)^2}-\frac{1}{(m-1)^2}&=&(\frac{1}{m-i}-\frac{1}{m-1})(\frac{1}{m-i}+\frac{1}{m-1})\\
&\le& \frac{2}{m-i}(\frac{1}{m-i}-\frac{1}{m-1})\\
&= &\frac{2(i-1)}{(m-i)^2(m-1)}\\
&\le& \frac{2(q-1)}{m-1}\frac1{(m-i)^2}\\
&\le& \frac{2(q-1)}{m-1}(\frac1{m-i-0.5}-\frac{1}{m-i+0.5}).
\end{eqnarray*}
Therefore, from \eqref{suff}
\begin{eqnarray*}
0\le \sum^q_{i=1}\frac{1}{(m-i)^2}-\frac{q}{(m-1)^2}&\le&\frac{2(q-1)}{m-1}(\frac1{m-q-0.5}-\frac{1}{m-0.5})\\
&=& \frac{2(q-1)q}{(m-1)(m-q-0.5)(m-0.5)}\\
&=&O(\frac{q^2}{m^2(m-q)})\\
&=&o(\xi(\frac{q}{m}))
\end{eqnarray*}
uniformly over $1\le 1<m$ as $m\to\infty$, proving \eqref{sum2}.

When $\delta=2$, we have from the above proof that
\[
\sum^q_{i=1}\frac{1}{(m-i)^2}\le \frac{q}{(m-1)^2}+O(\frac{q^2}{m^2(m-q)^2})=O(\frac{q}{(m-q)m}),
\]
that is,  \eqref{sum-delta} holds with $\delta=2$.   When $\delta>2$,
\[
\sum^q_{i=1}\frac{1}{(m-i)^{\delta}}\le \frac{1}{(m-q)^{\delta-2}}\sum^q_{i=1}\frac{1}{(m-i)^2}=O(\frac{1}{(m-q)^{\delta-2}}\frac{q}{(m-q)m})=O(\frac{q}{(m-q)^{\delta-1}m}),
\]
proving \eqref{sum-delta}.   \hfill$\blacksquare$

\begin{lemma}\label{lowerbound}
 Assume $p_n\ge 1$ and $k_n\ge 2$ are two sequences of integers such that $1\le p_n<\min_{1\le i\le k_n}n_i$ for all large $n$,  and $\min_{1\le i\le k_n}n_i\to\infty$ as $n\to\infty$. Then
\begin{equation}\label{liminf}
\liminf_{n\to\infty}\frac{\sigma_n^2}{p_n(p_n+1)(k_n-1)}\ge 1,
\end{equation}
and  for all large $n$
\begin{equation}\label{etan}
\eta(\frac{p_n}{n-k_n+1})\le \eta(\frac1{k_n})\le \eta(\frac12),
\end{equation}
where $\sigma_n^2$ and $\eta(x)$ are defined in \eqref{sigman} and \eqref{eta}, respectively.
\end{lemma}

\noindent{\it Proof.}  Without loss of generality, assume $n_1=\max_{1\le i\le k_n} n_i$. Since $\min_{1\le i\le k_n}n_i>p_n\ge 1$ and $k_n\ge 2$, we have
\[
n-k_n+1=\sum^{k_n}_{i=1}(n_i-1)+1\ge k_n(n_1-1)+1>n_1
\]
Also note that $\delta_n:=2\max_{1\le i\le k_n}(1-\frac{(n_i-1)^2}{n_i^2})=2(1-(1-\frac{1}{\min_{1\le i\le k_n}n_i})^2)\to 0$ as $n\to\infty$. From \eqref{eta-property}, $\eta(\frac{p_n}{n_i})\ge \eta(\frac{p_n}{n_1})\ge 1$ and $\eta(\frac{p_n}{n-k_n+1})\le \eta(\frac{p_n}{n_1})$.
We have for all large $n$
\begin{eqnarray*}
\sigma_n^2&=&p_n^2\Big(\sum^{k_n}_{i=1}\frac{(n_i-1)^2}{n_i^2}\eta(\frac{p_n}{n_i})-\frac{(n-k_n)^2}{(n-k_n+1)^2}\eta(\frac{p_n}{n-k_n+1})\Big)+p_n(k_n-1)\\
&\ge&p_n^2\Big(\sum^{k_n}_{i=1}\frac{(n_i-1)^2}{n_i^2}\eta(\frac{p_n}{n_1})-\eta(\frac{p_n}{n_1})\Big)+p_n(k_n-1)\\
&\ge&p_n^2\Big(\sum^{k_n}_{i=1}\frac{(n_i-1)^2}{n_i^2}-1\Big)\eta(\frac{p_n}{n_1})+p_n(k_n-1)\\
&\ge&p_n^2\Big(\sum^{k_n}_{i=1}(1-\frac{\delta_n}{2})-1\Big)+p_n(k_n-1)\\
&\ge&p_n^2(k_n-1)(1-\delta_n)+p_n(k_n-1)\\
&\ge&p_n(p_n+1)(k_n-1)(1-\delta_n),
\end{eqnarray*}
which implies \eqref{liminf}.

Since $n-k_n=\sum_{1\le i\le k_n}(n_i-1)\ge k_np_n$, we have
\begin{equation}\label{fraction}
\frac{p_n}{n-k_n+1}<\frac{1}{k_n}\le \frac12,
\end{equation}
which yields \eqref{etan} from the monotonicity of $\eta$ given in \eqref{eta-property}.
\hfill$\blacksquare$

\begin{lemma}\label{limitedcase} Assume there exists a constant $\delta\in (0,1)$ such that
 \begin{equation}\label{small}
 p_n\le \delta \min_{1\le i\le k_n}n_i~\mbox{for all large $n$,}~ \lim_{n\to\infty}\frac{\sqrt{k_n}}{\min_{1\le i\le k_n}n_i}=0,~\mbox{ and } \lim_{n\to\infty}\min_{1\le i\le k_n}n_i=\infty.
\end{equation}
Then
\begin{equation}\label{twobirds}
\lim_{n\to\infty}\frac{\mu_n-\bar\mu_n}{\sigma_n}=0,
\end{equation}
where $\mu_n$ and $\bar\mu_n$ are defined in \eqref{mun} and \eqref{barmun}, respectively.
\end{lemma}

\noindent{\it Proof.} Define
\[
g(x)=\psi(x)-\log x+\frac{1}{2x}, ~~~x>0.
\]
Then we have for any positive integers $m, q$ with $m>q\ge 1$
\begin{equation}\label{taumq}
\uptau (m, q):=\sum^q_{i=1}g(\frac{m-i}{2})=\sum^q_{i=1}\psi(\frac{m-i}{2})-\sum^q_{i=1}\log(\frac{m-i}{2})+\sum^q_{i=1}\frac{1}{m-i}.
\end{equation}



Since $n-k_n=\sum^{k_n}_{i=1}(n_i-1)\ge k_n(\min_{1\le i\le k_n}n_i-1)\to\infty$, we have from \eqref{fraction} that $n-k_n+1-p_n\ge \frac12(n-k_n)$ for all large $n$ and $n-k_n-p_n\to\infty$ as $n\to\infty$.

From \eqref{taumq}, \eqref{gamma0}, condition \eqref{small} and
Lemma~\ref{lowerbound} we have
\begin{eqnarray*} 
(n-k_n)\uptau (n-k_n+1, p_n)&=&(n-k_n)\sum^{p_n}_{j=1}g(\frac{n-k_n+1-j}{2})\nonumber\\
&=&O((n-k_n)\sum^{p_n}_{j=1}\frac{1}{(n-k_n+1-j)^2})\nonumber\\
&=&O(\frac{p_n}{n-k_n})\nonumber\\
&=&o(\sigma_n),
\end{eqnarray*}
and
\begin{eqnarray*} 
\sum^{k_n}_{i=1}(n_i-1)\uptau (n_i, p_n)&=&\sum^{k_n}_{i=1}(n_i-1)\sum^{p_n}_{j=1}g(\frac{n_i-j}{2})\nonumber\\
&=&\sum^{k_n}_{i=1}(n_i-1)\sum^{p_n}_{j=1}O(\frac{1}{(n_i-j)^2})\nonumber\\
&=&O(p_n)\sum^{k_n}_{i=1}\frac{1}{n_i}\nonumber\\
&=&O(\frac{p_nk_n}{\min_{1\le i\le k_n}n_i})\\
&=&o(p_n\sqrt{k_n})\\
&=&o(\sigma_n).
\end{eqnarray*}
Hence we have
\[
 (n-k_n)\uptau (n-k_n+1, p_n)-\sum^{k_n}_{i=1}(n_i-1)\uptau (n_i, p_n)=o(\sigma_n).
\]
In view of \eqref{mun}, \eqref{taumq}, and equations \eqref{mean} and \eqref{var} in Lemma~\ref{approximation12} we have
\begin{eqnarray*}
\mu_n&=&(n-k_n)\sum^{p_n}_{j=1}\log(\frac{n-k_n+1-j}{2})-\sum_{i=1}^{k_n}(n_i-1)\sum^{n_i}_{j=1}\log(\frac{n_i-j}{2})\\
&&-(n-k_n)\sum^{p_n}_{j=1}\frac{1}{n-k_n+1-j}+\sum_{i=1}^{k_n}(n_i-1)\sum^{p_n}_{j=1}\frac1{n_i-j}\\
&&+ (n-k_n)\uptau (n-k_n+1, p_n)-\sum^{k_n}_{i=1}(n_i-1)\uptau (n_i, p_n)\\
 &&+p_n\sum^{k_n}_{i=1}(n_i-1)\log(n_i-1)-p_n(n-k_n)\log(n-k_n)\\
 &=&(n-k_n)\sum^{p_n}_{j=1}\big(\log(n-k_n+1-j)-\log(n-k_n)\big)\\
 & &~-\sum_{i=1}^{k_n}(n_i-1)\sum^{p_n}_{j=1}\big(\log(n_i-j)-\log(n_i-1)\big)\\
&&~-(n-k_n)\sum^{p_n}_{j=1}\big(\frac{1}{n-k_n+1-j}-\frac{1}{n-k_n}\big)+\sum_{i=1}^{k_n}(n_i-1)\sum^{n_i}_{j=1}\big(\frac1{n_i-j}-\frac1{n_i-1}\big)\\
&&~+p_n(k_n-1) o(\sigma_n)\\
&=&(n-k_n)(p_n-n+k_n-\frac12)\log(1-\frac{p_n}{n-k_n+1})-\frac{(n-k_n)^2p_n}{n-k_n+1}\\
&&-\sum^{k_n}_{i=1}(n_i-1)(p_n-n_i+\frac12)\log(1-\frac{p_n}{n_i})+\sum^{k_n}_{i=1}\frac{(n_i-1)^2p_n}{n_i}\\
&&-(n-k_n)(-\frac{p_n}{n-k_n+1}-\log(1-\frac{p_n}{n-k_n+1}))+\sum^{k_n}_{j=1}(n_i-1)(-\frac{p_n}{n_i}-\log(1-\frac{p_n}{n_i}))\\
&&+p_n(k_n-1)+O(\frac{p_n}{n})+O(\sum^{k_n}_{i=1}\frac{p_n}{n_i})+o(\sigma_n)\\
&=&(n-k_n)(p_n-n+k_n+\frac12)\log(1-\frac{p_n}{n-k_n+1})-\sum^{k_n}_{i=1}(n_i-1)(p_n-n_i+\frac32)\log(1-\frac{p_n}{n_i})\\
&&+\sum^{k_n}_{i=1}\frac{(n_i-1)^2p_n}{n_i}-\frac{(n-k_n)^2p_n}{n-k_n+1}
+\frac{(n-k_n)p_n}{n-k_n+1}-\sum^{k_n}_{i=1}\frac{(n_i-1)p_n}{n_i}+p_n(k_n-1)+o(\sigma_n).
\end{eqnarray*}
Moreover, under condition \eqref{small}, we have
\begin{eqnarray*}
&&\sum^{k_n}_{i=1}\frac{(n_i-1)^2p_n}{n_i}-\frac{(n-k_n)^2p_n}{n-k_n+1}\\
&=&p_n\sum^{k_n}_{i=1}\frac{(n_i-1)^2-n_i^2}{n_i}-p_n\frac{(n-k_n)^2-(n-k_n+1)^2}{n-k_n+1} +\sum^{k_n}_{i=1}n_ip_n-(n-k_n+1)p_n\\
&=&p_n\sum^{k_n}_{i=1}(\frac1{n_i}-2)-p_n(\frac1{n-k_n+1}-2)+p_n(k_n-1)\\
&=&p_n\sum^{k_n}_{i=1}\frac1{n_i}-\frac{p_n}{n-k_n+1}-p_n(k_n-1)\\
&=&-p_n(k_n-1)+o(\sigma_n)
\end{eqnarray*}
and
\[
\frac{(n-k_n)p_n}{n-k_n+1}-\sum^{k_n}_{i=1}\frac{(n_i-1)p_n}{n_i}=-p_n(k_n-1)+\sum^{k_n}_{i=1}\frac{p_n}{n_i}-\frac{p_n}{n-k_n+1}=-p_n(k_n-1)+o(\sigma_n).
\]
Then we conclude \eqref{twobirds}. \hfill$\blacksquare$


Let $f$ be any function defined over $(0,\infty)$. For integers $q$
and $m$ with $1\le q<m$, define the function
\begin{equation}\label{Sfunction}
\Lambda_{f,\; m, q}(t)=\sum^q_{i=1}f(\frac{m-i}{2}+t),
~~t>-\frac{m-q}{2}.
\end{equation}
We need this definition in the following proofs.


\begin{lemma}\label{G} For any integers $1\le q< m$, define
\[
G_{m,q}(x)=\sum^q_{i=1}\Big(s(\frac{(m-i)}{2}+x)-s(\frac{m-1}{2}+x)\Big), ~~~x>\frac{m-q}{2},
\]
where the function $s(x)$ is defined in Lemma~\ref{loggamma}. Then
\begin{eqnarray}\label{GG}
&&G_{m,q}(\frac{(m-1)t}{2})-G_{m,q}(0)\\
&=&\Big(\Lambda_{s',\;m,q}(0)-qs'(\frac{m-1}{2})\Big)\frac{(m-1)t}{2}+\Big(\xi(\frac{q}{m})(1+o(1))-\frac{q}{(m-1)^2}\Big)\frac{(m-1)^2t^2}{8}\nonumber
\end{eqnarray}
and
\begin{eqnarray}\label{ss}
&&\Lambda_{s,\;m,q}(\frac{(m-1)t}{2})-\Lambda_{s,\;m,q}(0)\\
&=&\Big(\Lambda_{s',\;m,q}(0)-qs'(\frac{m-1}{2})\Big)\frac{(m-1)t}{2}+\Big(\xi(\frac{q}{m})(1+o(1))-\frac{q}{(m-1)^2}\Big)\frac{(m-1)^2t^2}{8}\nonumber\\
& &~~~+\frac{(m-1)q}{2}(1+t)\log(1+t)+\frac{(m-1)q}{2}(1+\log\frac{m-1}{2})t-\frac{q}{2}\log(1+t)\nonumber
\end{eqnarray}
hold uniformly over $|t|\le
(\varepsilon+\frac{1}{m-q})\frac{m-q}{2(m-1)}$, $1\le q<m$ as
$\varepsilon\to 0$ and  $m\to\infty$.
\end{lemma}

\noindent{\it Proof.} We will apply Taylor's theorem to expand a smooth function $f(x)$, that is,
that is, for any $x>0$, $x+y>0$ there exists a constant $c$ between $x$ and $x+y$ such that
\begin{equation}\label{taylor}
f(x+y)=f(x)+f'(x)y+f''(x)\frac{y^2}{2}+f'''(c)\frac{y^3}{6}.
\end{equation}
In our application here, we will assume $x\ge \frac{1}{2}$, and $|y|\le \frac{1}{2}x$.

We note that
\begin{equation}\label{hh}
\frac{m-i}{2}\ge \frac14, ~~~
\frac12\le \frac{\frac{(m-i)}{2}+ \frac{(m-1)t}{2}}{\frac{(m-i)}{2}} \le \frac32,  \
\end{equation}
when $|t|\le \frac{m-q}{2(m-1)}$,  $1\le i\le q<m$.

Now it follows from \eqref{taylor} that for each $1\le i\le q$, there a constant $c_i$ between $\frac{m-i}{2}$ and $\frac{(m-i)}{2}+\frac{(m-1)t}{2}$ such that
\begin{eqnarray*}
g_i(t):&=&s(\frac{(m-i)}{2}+\frac{(m-1)t}{2})-s(\frac{m-1}{2}+\frac{(m-1)t}{2})-(s(\frac{m-i}{2})-s(\frac{m-1}{2}))\\
&=&\big(s'(\frac{m-i}{2})-s'(\frac{m-1}{2})\big)\frac{(m-1)t}{2}+\big(s''(\frac{m-i}{2})-s''(\frac{m-1}{2})\big)\frac{(m-1)^2t^2}{8}\\
&&~~~~+\big(s'''(\frac{(m-i)}{2}+c_i)-s'''(\frac{m-1}{2}+c_i)\big)\frac{(m-1)^3t^3}{48}.
\end{eqnarray*}
Set
\[
R_i=\big(s'''(\frac{(m-i)}{2}+c_i)-s'''(\frac{m-1}{2}+c_i)\big)\frac{(m-1)^3t^3}{48}.
\]
Then
\begin{eqnarray*}
G_{m,q}(\frac{(m-1)t}{2})&=&\sum^q_{i=1}(s(\frac{m-i}{2})-s(\frac{m-1}{2}))+\sum^q_{i=1}g_i(t)\\
&=&G_{m,q}(0)+\sum^q_{i=1}g_i(t)\\
&=&G_{m,q}(0)+\Big(\Lambda_{s',\;m,q}(0)-qs'(\frac{m-1}{2})\Big)\frac{(m-1)t}{2}\\
&&+\Big(\Lambda_{s'',\;m,q}(0)-qs''(\frac{m-1}{2})\Big)\frac{(m-1)^2t^2}{8}+\sum^q_{i=1}R_i.
\end{eqnarray*}

From \eqref{s-diff}, \eqref{hh}, and Lemma~\ref{x-delta} we have
\begin{eqnarray*}
|R_i|&\le&\left(\frac{1}{(\frac{m-i}{2}+c_i)^2}-\frac{1}{(\frac{m-1}{2}+c_i)^2}
+\frac{1}{(\frac{m-i}{2}+c_i)^3}-\frac{1}{(\frac{m-1}{2}+c_i)^3}\right)\frac{(m-1)^3|t|^3}{48}\\
&\le& \left(\frac{2i}{(\frac{m-i}{2}+c_i)^2(\frac{m-1}{2}+c_i))}+\frac{3i}{(\frac{m-i}{2}+c_i)^3(\frac{m-1}{2}+c_i))}\right)\frac{(m-1)^3|t|^3}{48}\\
&\le&\frac{64q}{(m-i)^2(m-1)}\frac{(m-1)^3|t|^3}{48}\\
&\le&\frac{2(m-1)^2q|t|^3}{(m-i)^2}
\end{eqnarray*}
and, thus from \eqref{sum-delta}, \eqref{suff}, \eqref{eta} and
\eqref{eta-property}
\begin{eqnarray*}
\sum^q_{i=1}|R_i|&\le& 2(m-1)^2q|t|^3\sum^q_{i=1}\frac{1}{(m-i)^2}\\
&=&O(\frac{(m-1)q^2|t|^3}{(m-q)}\\
&=&O((\varepsilon+\frac{1}{m-q})\frac{m-q}{2(m-1)}\frac{(m-1)q^2|t|^2}{(m-q)})\\
&=&O( \frac{\varepsilon q^2}{(m-1)^2}  +\frac{q^2}{(m-1)^2(m-q)})(m-1)^2t^2\\
&=&o(\xi(\frac{q}{m}))(m-1)^2t^2
\end{eqnarray*}
holds uniformly over $|t|\le
(\varepsilon+\frac{1}{m-q})\frac{m-q}{2(m-1)}$, $1\le q<m$ as
$\varepsilon\to 0$ and  $m\to\infty$. Furthermore, in view of
\eqref{s-diff}, \eqref{m-i-m} and \eqref{sum2} we have
\begin{eqnarray*}
\Lambda_{s'',\;m,q}(0)-qs''(\frac{m-1}{2})&=&2\sum^q_{i=1}\Big(\frac{1}{m-i}-\frac{1}{m-1}\Big)+2\sum^q_{i=1}\Big(\frac{1}{(m-i)^2}-\frac{1}{(m-1)^2}\Big)\\
&=&\xi(\frac{q}{m})-\frac{q}{(m-1)^2}+o(\xi(\frac{q}{m}))\\
\end{eqnarray*}
uniformly over $1\le q< m$ as $m\to\infty$.    \eqref{GG} can be easily obtained by combining all results above.

Review the definition of $s(x)$ in Lemma~\ref{loggamma} and the definition of $\Lambda$ in \eqref{Sfunction}. Since
 $\Lambda_{s,\;m,q}(x)=G_{m,q}(x)+qs(\frac{m-1}{2}+x)$, we have
 \[
 \Lambda_{s,\;m,q}(\frac{(m-1)t}{2})-\Lambda_{s, m,q}(0)=G_{m,q}(\frac{(m-1)t}{2})-G_{m,q}(0)+q\big(s(\frac{(m-1)(1+t)}{2})-s(\frac{m-1}{2})\big).
 \]
By using \eqref{GG} and the fact that
\begin{eqnarray*}
&&s(\frac{(m-1)(1+t)}{2})-s(\frac{m-1}{2})\\
&=&(\frac{(m-1)(1+t)}{2}-\frac{1}{2})\log(\frac{(m-1)(1+t)}{2})-(\frac{m-1}{2}-\frac{1}{2})\log(\frac{m-1}{2})+\frac{(m-1)t}{2}\\
&=&\frac{m-1}{2}(1+t)\log(1+t)+\frac{m-1}{2}(1+\log\frac{m-1}{2})t-\frac{1}{2}\log(1+t),
\end{eqnarray*}
we conclude \eqref{ss}.
\hfill$\blacksquare$

\begin{lemma}\label{gamma-Log} We have
\begin{eqnarray*}
& &\log\frac{\Gamma_q(\frac{m-1}{2}(1+t))}{\Gamma_q(\frac{m-1}2)}\\
&=&\Big(\Lambda_{\psi,\;m,q}(0)-qs'(\frac{m-1}{2})\Big)\frac{(m-1)t}{2}+\Big(\xi(\frac{q}{m})(1+o(1))-\frac{q}{(m-1)^2}\Big)\frac{(m-1)^2t^2}{8}\nonumber\\
& &~~~+\frac{(m-1)q}{2}(1+t)\log(1+t)+\frac{(m-1)q}{2}(1+\log\frac{m-1}{2})t-\frac{q}{2}\log(1+t)
\end{eqnarray*}
uniformly over $|t|\le
(\varepsilon+\frac{1}{m-q})\frac{m-q}{2(m-1)}$, $1\le q<m$ as
$\varepsilon\to 0$ and  $m\to\infty$.
\end{lemma}

\noindent {\it Proof.} By the definition of $\Gamma_q$ function in \dref{gammap}, we have
\begin{equation}\label{logGamma}
\log\frac{\Gamma_q(\frac{m-1}{2}(1+t))}{\Gamma_p(\frac{m-1}2)}=\sum_{i=1}^q\log\frac{\Gamma(\frac{m-i}{2}+\frac{(m-1)t}{2})}{\Gamma (\frac{m-i}{2})}=\Lambda_{\log\Gamma,\; m, q}(\frac{(m-1)t}{2})-\Lambda_{\log\Gamma,\; m, q}(0).
\end{equation}
The function is well defined when $t>-\frac{m-q}{m-1}$.

 For function $h(x)$ defined in \eqref{h(x)}, we apply Taylor's formula with remainder
 \begin{equation}\label{hhi}
 h(\frac{m-i}{2}+\frac{(m-1)t}{2})-h(\frac{m-i}{2})=h'(\frac{m-i}{2})\frac{(m-1)t}{2}+h''(c_i)\frac{(m-1)^2t^2}{8}
 \end{equation}
 where $c_i$ is a number between $\frac{m-i}{2}+\frac{(m-1)t}{2}$ and $\frac{m-i}{2}$. Using Lemma~\ref{loggamma} we have for some $C>0$
\[
|h''(c_i)|\le \frac{C}{(m-i)^3},~~~1\le i\le q<m<\infty.
\]
Summing up over $1\le i\le q$ on both sides of \eqref{hhi} yields
\begin{eqnarray}\label{Delta-h}
\Lambda_{h,\;m,q}(\frac{(m-1)t}{2})-\Lambda_{h,\;m,q}(0)&=&\Lambda_{h',\;m,q}(0)\frac{(m-1)t}{2}+O(\sum^q_{i=1}\frac{1}{(m-i)^3})(m-1)^2t^2\nonumber\\
&=&\Lambda_{h',\;m,q}(0)\frac{(m-1)t}{2}+O(\frac{q}{(m-q)^2m})(m-1)^2t^2\nonumber\\
&=&\Lambda_{h',\;m,q}(0)\frac{(m-1)t}{2}+o(\xi(\frac{q}{m}))(m-1)^2t^2
\end{eqnarray}
uniformly over $|t|\le
(\varepsilon+\frac{1}{m-q})\frac{m-q}{2(m-1)}$, $1\le q<m$ as
$\varepsilon\to 0$ and  $m\to\infty$. In the above estimation, we
have used \eqref{sum-delta} from Lemma~\ref{approximation12} and
\eqref{suff2} from Lemma~\ref{oxi}.

From \eqref{h(x)}, we have
\[
\Lambda_{\log\Gamma,\;m,q}(\frac{(m-1)t}{2})=\Lambda_{h,\;m,q}(\frac{(m-1)t}{2})+\Lambda_{s,\;m,q}(\frac{(m-1)t}{2}).
\]

%


Now by combining \eqref{logGamma}, \eqref{Delta-h}, \eqref{ss} and the above equations we have
\begin{eqnarray*}
&&\log\frac{\Gamma_q(\frac{m-1}{2}(1+t))}{\Gamma_q(\frac{m-1}2)}\\
&=&\Lambda_{h,\;m,q}(\frac{(m-1)t}{2})-\Lambda_{h,\;m,q}(0)+\Lambda_{s,\;m,q}(\frac{(m-1)t}{2})-\Lambda_{s,\;m,q}(0)\\
&=&\Lambda_{h',\;m,q}(0)\frac{(m-1)t}{2}+o(\xi(\frac{q}{m}))(m-1)^2t^2\\
& &+\Big(\Lambda_{s',\;m,q}(0)-qs'(\frac{m-1}{2})\Big)\frac{(m-1)t}{2}+\Big(\xi(\frac{q}{m})(1+o(1))-\frac{q}{(m-1)^2}\Big)\frac{(m-1)^2t^2}{8}\nonumber\\
& &~~~+\frac{(m-1)q}{2}(1+t)\log(1+t)+\frac{(m-1)q}{2}(1+\log\frac{m-1}{2})t-\frac{q}{2}\log(1+t)\\
&=&\Big(\Lambda_{h',\;m,q}(0)+\Lambda_{s',\;m,q}(0)-qs'(\frac{m-1}{2})\Big)\frac{(m-1)t}{2}+\Big(\xi(\frac{q}{m})(1+o(1))-\frac{q}{(m-1)^2}\Big)\frac{(m-1)^2t^2}{8}\nonumber\\
& &~~~+\frac{(m-1)q}{2}(1+t)\log(1+t)+\frac{(m-1)q}{2}(1+\log\frac{m-1}{2})t-\frac{q}{2}\log(1+t)\\
&=&\Big(\Lambda_{\psi,\;m,q}(0)-qs'(\frac{m-1}{2})\Big)\frac{(m-1)t}{2}+\Big(\xi(\frac{q}{m})(1+o(1))-\frac{q}{(m-1)^2}\Big)\frac{(m-1)^2t^2}{8}\nonumber\\
& &~~~+\frac{(m-1)q}{2}(1+t)\log(1+t)+\frac{(m-1)q}{2}(1+\log\frac{m-1}{2})t-\frac{q}{2}\log(1+t)
\end{eqnarray*}
uniformly over $|t|\le
(\varepsilon+\frac{1}{m-q})\frac{m-q}{2(m-1)}$, $1\le q<m$ as
$\varepsilon\to 0$ and  $m\to\infty$.  This completes the proof of
the lemma. \hfill$\blacksquare$

Finally, we state the following lemma to express the moments of $W_n$ as a product of $\Gamma_p$ functions, which is available on page 302 in Muirhead~\cite{muirhead82}.

\begin{lemma}\label{moment}
Let $\Lambda_n^*$ be as in \dref{Bartlett}. Set
\begin{equation}\label{wn}
W_n=\frac{\prod_{i=1}^k|A|^{(n_i-1)/2}}{|A|^{(n-k)/2}}=\Lambda_n^*\cdot\frac{\prod_{i=1}^k(n_i-1)^{(n_i-1)p/2}}{(n-k)^{(n-k)p/2}}.
\end{equation}
Assume $n_i>p$ for $1\le i \le k$. Then under $H_0$ in \dref{test}, we have
\[
E(W_n^t)=\frac{\Gamma_p\left({1 \over 2}(n-k)\right)}{\Gamma_p\left({1 \over 2}(n-k)(1+t)\right)} \cdot \prod_{i=1}^{k}\frac{\Gamma_p\left({1 \over 2}(n_i-1)(1+t)\right)}{\Gamma_p\left({1 \over 2}(n_i-1)\right)},
\]
for all $t > \max_{1\le i \le k}{p-1 \over n_i-1}-1$, where $\Gamma_p$ is defined as in \dref{gammap}.
\end{lemma}

\subsection{Proofs of the main results}\label{proofofmain}

\noindent {\bf Proof of Theorem \ref{thm1}.}

We first show the sufficiency, that is,  \eqref{CLT} holds under
condition $\max(p_n, k_n)\to\infty$.

 By \dref{wn}, we write
\[
-2\log \Lambda_n^*=-2\log
W_n+p_n\Big(\sum_{i=1}^{k_n}(n_i-1)\log(n_i-1)-(n-k_n)\log(n-k_n)\Big).
\]

Let
\[
\mu_n'=\mu_n-p_n\Big(\sum_{i=1}^{k_n}(n_i-1)\log(n_i-1)-(n-k_n)\log(n-k_n)\Big).
\]
that
\[
\frac{-2\log W_n-\mu_n'}{\sigma_n} \stackrel{d}{\to} N(0,1),
\]
as $n \rightarrow \infty$. We will use the moment generating function of $-2\log W_n$ and all we need is to prove that
\[
E\exp\left\{\frac{-2\log W_n-\mu_n'}{\sigma_n}w\right\}
=E\left(W_n^{-2w/\sigma_n}\right)\exp(\mu_n'(-w/\sigma_n))\rightarrow
e^{w^2/2}
\]
for all $|w|\le 1/8$. Now we write $t_n=-2w/\sigma_n$. We need to show
\[
\lim_{n\to\infty}\left(\log E\left(W_n^{t_n}\right)+\frac12\mu_n't_n\right)=\frac{w^2}2.
\]

When $|w|\le \frac18$, $|t_n|\le \frac{1}{4\sigma_n}$.  According to Lemma~\ref{lowerbound}, for all large $n$,  $\sigma_n\ge \frac12p_n\sqrt{k_n-1}$ and thus $|t_n|\le \frac{1}{2p_n\sqrt{k_n-1}}$. In this case,
\[
t_n>-\frac{1}{2p_n\sqrt{k_n-1}}>-\frac{1}{p_n}>\frac{p_n-1}{p_n}-1\ge \frac{p_n-1}{\min_{1\le i\le k_n}n_i}-1=\max_{1\le i\le k_n}\frac{p_n-1}{n_i}-1.
\]
Therefore, we can apply Lemma~\ref{moment} for all large $n$ to get
\begin{equation}\label{logewn}
\log E(W_n^{t_n})=\sum_{i=1}^{k_n}\log \frac{\Gamma_{p_n}\left({1 \over 2}(n_i-1)(1+t_n)\right)}{\Gamma_{p_n}\left({1 \over 2}(n_i-1)\right)}
-\log \frac{\Gamma_{p_n}\left({1 \over 2}(n-k_n)(1+t_n)\right)}{\Gamma_{p_n}\left({1 \over 2}(n-k_n)\right)}.
\end{equation}
The following equations will be used in approximating $\log E(W_n^{t_n})$ in \eqref{logewn}:
\[
n-k_n=\sum^{k_n}_{i=1}(n_i-1),
\]
\begin{equation}\label{sumofs'}
\sum^{k_n}_{i=1}(n_i-1)s'(\frac{n_i-1}{2})-(n-k_n)s'(\frac{n-k_n}{2})=\sum^{k_n}_{i=1}(n_i-1)\log(n_i-1)-(n-k_n)\log(n-k_n)-(k_n-1),
\end{equation}
\[
\frac{(n-k_n)^2p_n^2}{(n-k_n+1)^2}\eta(\frac{p_n}{n-k_n+1})=O(\sigma_n^2)
\]
from Lemma~\ref{lowerbound},
\[
\sum^{k_n}_{i=1}(n_i-1)^2\xi(\frac{p_n}{n_i})\le \sigma^2_n-\frac{(n-k_n)^2p_n^2}{(n-k_n+1)^2}\eta(\frac{p_n}{n-k_n+1})=O(\sigma_n^2).
\]
and  $t_n\to 0$ as $n\to\infty$ since $\max(k_n,p_n)\to\infty$.


For any $m>p_n$,
\begin{eqnarray*}
|t_n|&\le&\frac{1}{2p_n\sqrt{k_n-1}}\frac{2(m-1)}{m-p_n} \frac{m-p_n}{2(m-1)}\\
&<&\frac{1}{\sqrt{k_n-1}}(\frac{m-p_n+p_n)}{(m-p_n}p_n)\frac{m-p_n}{2(m-1)}\\
&=&\frac{1}{\sqrt{k_n-1}}\big(\frac1{p_n}+\frac1{m-p_n}\big)\frac{m-p_n}{2(m-1)}\\
&=&\big(\frac1{p_n\sqrt{k_n-1}}+\frac1{m-p_n}\big)\frac{m-p_n}{2(m-1)}.
\end{eqnarray*}
Set $\varepsilon_n=\frac{1}{p_n\sqrt{k_n-1}}$. Then
$\lim_{n\to\infty}\varepsilon_n=0$ since $\max(p_n,k_n)\to\infty$.
Clearly, we have
\[
|t_n|\le \big(\varepsilon_n+\frac1{m-p_n}\big)\frac{m-p_n}{2(m-1)}
\]
for $m=n_i$, $1\le i\le k_n$ and  $m=n-k_n+1$.  This implies
Lemma~\ref{gamma-Log} holds uniformly over $m=n_i$, $1\le i\le k_n$
and $m=n-k_n$ with $q=p_n$ as $n\to\infty$. Therefore, from
\eqref{logewn}
\begin{eqnarray*}
&&\log E(W_n^{t_n})\\
&=&\sum_{i=1}^{k_n}\Big(\Lambda_{\psi,\;n_i,p_n}(0)-p_ns'(\frac{n_i-1}{2})\Big)\frac{(n_i-1)t_n}{2}-\Big(\Lambda_{\psi,\;n-k_n+1,p_n}(0)-p_ns'(\frac{n-k_n}{2})\Big)\frac{(n-k_n)t_n}{2}\\
&&~~~+\sum_{i=1}^{k_n}\Big(\xi(\frac{p_n}{n_i})(1+o(1))-\frac{p_n}{(n_i-1)^2}\Big)\frac{(n_i-1)^2t_n^2}{8}\\
&&~~~-\Big(\xi(\frac{p_n}{n-k_n+1})(1+o(1))-\frac{p_n}{(n-k_n)^2}\Big)\frac{(n-k_n)^2t_n^2}{8}\\
&&~~~+\sum_{i=1}^{k_n}\Big(\frac{(n_i-1)p_n}{2}(1+t_n)\log(1+t_n)+\frac{(n_i-1)p_n}{2}(1+\log\frac{n_i-1}{2})t_n-\frac{p_n}{2}\log(1+t_n)\Big)\\
&&~~~-\Big(\frac{(n-k_n)p_n}{2}(1+t_n)\log(1+t_n)+\frac{(n-k_n)p_n}{2}(1+\log\frac{n-k_n}{2})t_n-\frac{p_n}{2}\log(1+t_n)\Big)\\
&=&\frac{t_n}{2}\Big\{\sum_{i=1}^{k_n}(n_i-1)\Lambda_{\psi,\;n_i,p_n}(0)- (n-k_n)\Lambda_{\psi,\;n-k_n+1,p_n}(0)\Big\}\\
&&~~~-\frac{p_nt_n}{2}\Big\{\sum^{k_n}_{i=1}(n_i-1)s'(\frac{n_i-1}{2})-(n-k_n)s'(\frac{n-k_n}{2})\big)\Big\}\\
&&~~~+\frac{t_n^2}{8}\Big\{\sum_{i=1}^{k_n}(n_i-1)^2\xi(\frac{p_n}{n_i})-(n-k_n)^2\xi(\frac{p_n}{n-k_n+1})-p_n(k_n-1)\Big\}(1+o(1))\\
&&~~~+\frac{p_nt_n}{2}\Big\{\sum_{i=1}^{k_n}(n_i-1)\log\frac{n_i-1}{2}-(n-k_n)\log\frac{n-k_n}{2}\Big\}\\
&&~~~-\frac{p_n(k_n-1)}{2}\log(1+t_n)\\
&=&: I_1+I_2+I_3+I_4+I_5.
\end{eqnarray*}

Since
\[
I_1=-\frac12\mu_n't_n,
\]
\[
I_2+I_4=\frac{1}{2}p_n(k_n-1)t_n
\]
from \eqref{sumofs'},
and
\[
I_5=-\frac{p_n(k_n-1)}{2}\log(1+t_n)=-\frac12p_n(k_n-1)t_n+p_n(k_n-1)\frac{t_n^2}{4}(1+o(1)).
\]
by Taylor's expansion, we have
\begin{eqnarray*}
&&\log E(W_n^{t_n})+\frac12\mu_n't_n\\
&=&\frac{t_n^2}{8}\Big(\sum_{i=1}^{k_n}(n_i-1)^2)\xi(\frac{p_n}{n_i})-(n-k_n)^2\xi(\frac{p_n}{n-k_n+1})+p_n(k_n-1)\Big)(1+o(1))\\
&=&\frac18\Big(\frac{-2w}{\sigma_n}\Big)^2\sigma_n^2(1+o(1))\\
&\to&\frac{w^2}{2}
\end{eqnarray*}
as $n\to\infty$.

To prove  the necessity, we need to show that \eqref{CLT} implies
$\lim_{n\to\infty}\max(p_n, k_n)=\infty$.  If condition
$\lim_{n\to\infty}\max(p_n, k_n)=\infty$ is not sure, then we see
that exists a subsequence of $n$, say, $\{n'\}$ such that $p_{n'}=p$
and $k_{n'}=k$ for all large $n'$, where $p$ and $k$ are two fixed
integers. Then it follows from the classic chi-square approximation
\eqref{chiapprox} that $-2\log\Lambda_{n'}^*$ converges in
distribution to a chi-square distribution with
$f=\frac12p(p+1)(k-1)$ degrees of freedom as $n'\to\infty$. This is
contradictory to \eqref{CLT}.

 This completes the proof of the theorem.   \hfill$\blacksquare$

\vspace{20pt}

\noindent{\it Proof of Theorem~\ref{thm-chisq}.}

To prove \eqref{chisquare}, it suffices to show that
for any subsequence $\{n'\}$ of $\{n\}$,
there is a further subsequence $\{n''\}$ such that \eqref{chisquare} holds along $\{n''\}$.  The subsequence
$\{n''\}$ can be selected in a way that both the limits of  $k_{n''}$ and $p_{n''}$  exist, and both the limits
can be infinity.    For the sake of simplicity,  we can assume both the limits of  $k_{n}$ and $p_{n}$  exist along
the entire sequence and prove \eqref{chisquare} holds. Note that if the limit of a sequence of integers is finite, the sequence takes a constant value ultimately. For this reason, it suffices to show \eqref{chisquare} under conditions  $\min_{1\le i\le k_n}n_i>p_n$ and $\min_{1\le i\le k_n}n_i\to\infty$, and each of the following two assumptions:
\newline
\noindent\textbf{Assumption 1:}  $\lim_{n\to\infty}\max(p_n,
k_n)=\infty$;
\newline
\noindent \textbf{Assumption 2:} $p_n=p$ and $k_n=k$ for all large
$n$, where both $p$ and $k$ are fixed integers.

Assume Assumption 1 holds. In this case, Theorem~\ref{thm1} holds, and $f_n=\frac12{p_n(p_n+1)}{k_n-1} \to\infty$ as $n\to\infty$. By using Theorem~\ref{thm1} and applying the central limit theorem to $\chi^2_{f_n}$ we have
\[
\sup_{-\infty<x<\infty}|P(\frac{-2\log\Lambda_n-\mu_n}{\sigma_n}\le x)-\Phi(x)|\to 0\mbox{ and } \sup_{-\infty<x\infty}|P(\frac{\chi^2_{f_n}-f_n}{\sqrt{2f_n}}\le x)-\Phi(x)|\to 0
\]
  as $n\to\infty$, where $\Phi(x)$ denotes the cumulative distribution function for the standard normal random variable. Therefore,
\begin{eqnarray*}
&&\sup_{-\infty<x<\infty}|P(Z_n\le x)-P(\chi^2_{f_n}\le x)|\\
&=&\sup_{-\infty<x<\infty}|P(\frac{Z_n-f_n}{\sqrt{2f_n}}\le x)-P(\frac{\chi^2_{f_n}-f_n}{\sqrt{2f_n}}\le x)|\\
&\le &\sup_{-\infty<x<\infty}|P(\frac{Z_n-f_n}{\sqrt{2f_n}}\le x)-\Phi(x)|+\sup_{-\infty<x<\infty}|P(\frac{\chi^2_{f_n}-f_n}{\sqrt{2f_n}}\le x)-\Phi(x)|\\
&\to& 0~~~\mbox{ as }n\to\infty,
\end{eqnarray*}
proving \eqref{chisquare}.

Under Assumption 2,  we have that $f_n=f=\frac{1}{2}p(p+1)(k-1)$, as defined in \eqref{fn}, is a fixed integer.  It suffices to show that $Z_n$ converges in distribution to $\chi^2_f$.  From \eqref{zn}, if we are able to show
\begin{equation}\label{poissonmean}
\frac{2f_n}{\sigma_n^2}\to 1~~~\mbox{ and }~~~f_n-\mu_n\sqrt{\frac{2f_n}{\sigma_n^2}}\to 0
\end{equation}
as $n\to\infty$, then $Z_n=(-2\log \Lambda_n^*)(1+o(1))+o(1)$, which concludes the desired result by using \eqref{chiapprox}.

It is easy to see that
\begin{eqnarray*}
\lim_{n\to\infty}\sigma_n^2&=&\lim_{n\to\infty}\Big(\sum^{k}_{i=1}(n_i-1)^2(\frac{p}{n_i})^2\eta(\frac{p}{n_i})
-(n-k_n)^2(\frac{p}{n-k_n})^2\eta(\frac{p}{n-k_n+1})\Big)+p(k-1)\\
&=&\sum^k_{i=1}p^2\eta(0)-p^2\eta(0)+p(k-1)\\
&=&p^2(k-1)+p(k-1)\\
&=&p(p+1)(k-1).
\end{eqnarray*}

Next, we need to estimate $\mu_n$. Since Assumption 2 implies
condition \eqref{small}, we can estimate $\bar\mu_n$ and then apply
\eqref{twobirds}.  By using Taylor's expansion for $\bar\mu_n$ given
in \eqref{barmun}, we can show
$\bar\mu_n=\frac{1}{2}p(p+1)(k-1)+o(1)$.  We omit the details here.
Therefore, we have
$\mu_n=\bar\mu_n+\sigma_n(\frac{\mu_n-\bar\mu}{\sigma_n})\to
\frac12p(p+1)(k-1)$ as $n\to\infty$.  Then \eqref{poissonmean}
follows easily from the limits of $\sigma_n^2$ and $\mu_n$.
\hfill$\blacksquare$

\vspace{10pt}

\noindent{\it Proof of Theorem~\ref{thm-special}.} Part of (b) of
the theorem follows from Theorem~\ref{thm1} and equation
\eqref{twobirds} in Lemma~\ref{limitedcase}. Part (a) follows from
the same lines in the proof of Theorem~\ref{thm-chisq} by using Part
(b) of Theorem~\ref{thm-special}.  We omit the details.
\hfill$\blacksquare$


\vspace{10pt}

\noindent\textbf{Supplementary Material}

The Supplementary Material contains Figures S2 and S3.

\vspace{10pt}

\noindent\textbf{Acknowledgements}

The authors would like to thank the Editor, Associate Editor, and
referees for reviewing the manuscript and providing valuable
comments.  The research of Yongcheng Qi was supported in part by NSF
Grant DMS-1916014.

\vspace{10pt}

\noindent\textbf{Conflict of Interest Statement}

On behalf of all authors, the corresponding author states that there
is no conflict of interest.


\baselineskip 12pt
\def\ref{\par\noindent\hangindent 25pt}

\newpage







\def\C{\mathbb{C}}


\setcounter{page}{1}

\thispagestyle{empty}

\noindent {\bf \large Asymptotic Distributions for Likelihood Ratio
Tests \\
for the Equality of Covariance Matrices}

\vspace{20pt}
 \noindent{\bf Wenchuan Guo$^a$, Yongcheng Qi$^b$}

\vspace{10pt}

\noindent $^a$Biometrics \& Data Sciences, Bristol Myers Squibb,
3551 Lawrenceville Princeton, Lawrence Township, NJ 08543, USA.

\vspace{10pt}

\noindent $^b$Department of Mathematics and Statistics, University
of Minnesota Duluth, 1117 University Drive, Duluth, MN 55812, USA.

\date{\today}

\vspace{20pt}




\vspace{0.3cm}

\vspace{.55cm} \centerline{\bf Supplementary Material}
\vspace{.55cm} \fontsize{9}{11.5pt plus.8pt minus .6pt}\selectfont

\noindent We present additional plots for histograms of $-2\rho
\log\Lambda_n^*$ (Chisq), $(-2\log\Lambda_n^*-\mu_n)/\sigma_n$
(CLT), and  $Z_n$ (ALRT).

\par

\setcounter{figure}{1} \setcounter{equation}{0}
\def\thefigure{S\arabic{figure}}

\fontsize{12}{14pt plus.8pt minus .6pt}\selectfont


\begin{figure}[H]
\centering
\includegraphics[height=2.0in,width=6in]{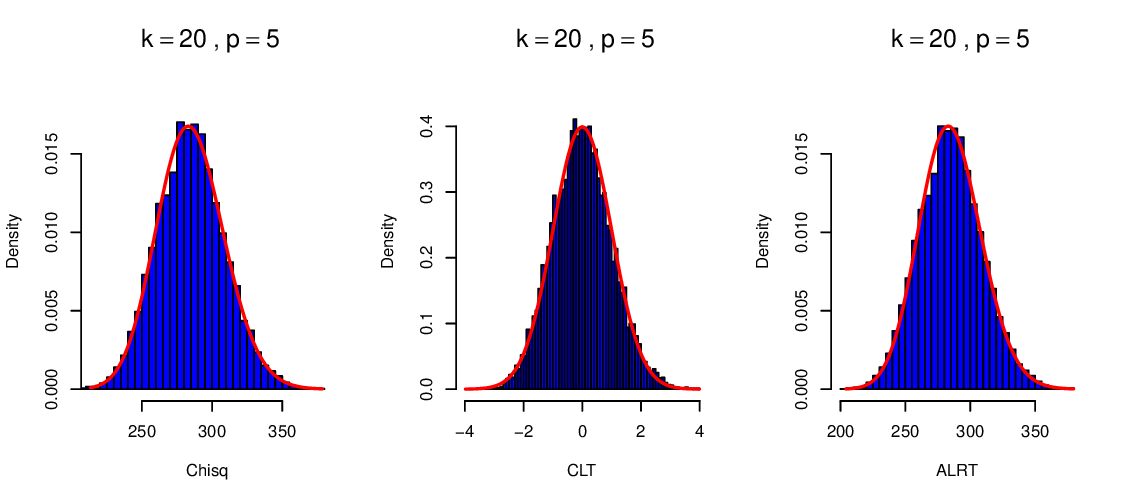}
\includegraphics[height=2.0in,width=6in]{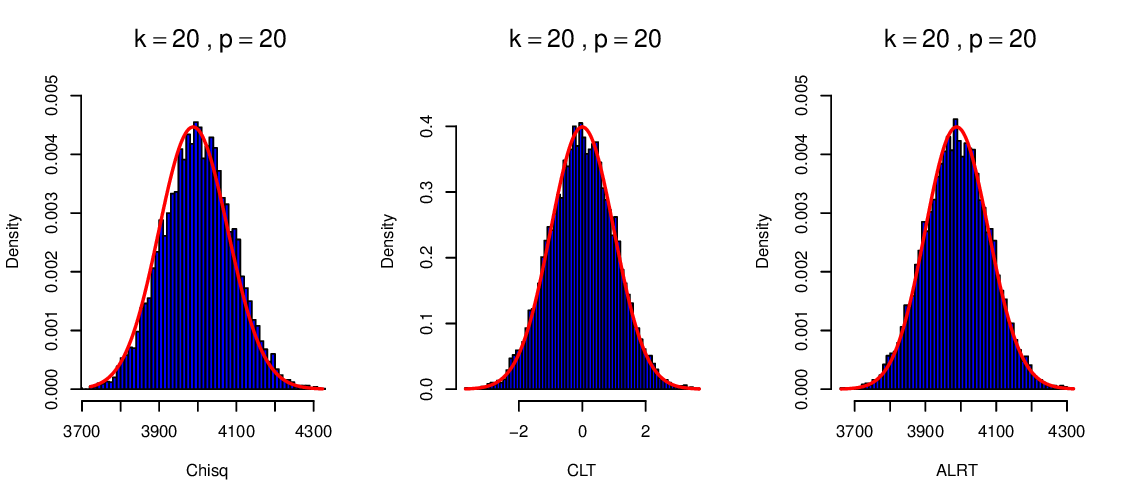}
\includegraphics[height=2.0in,width=6in]{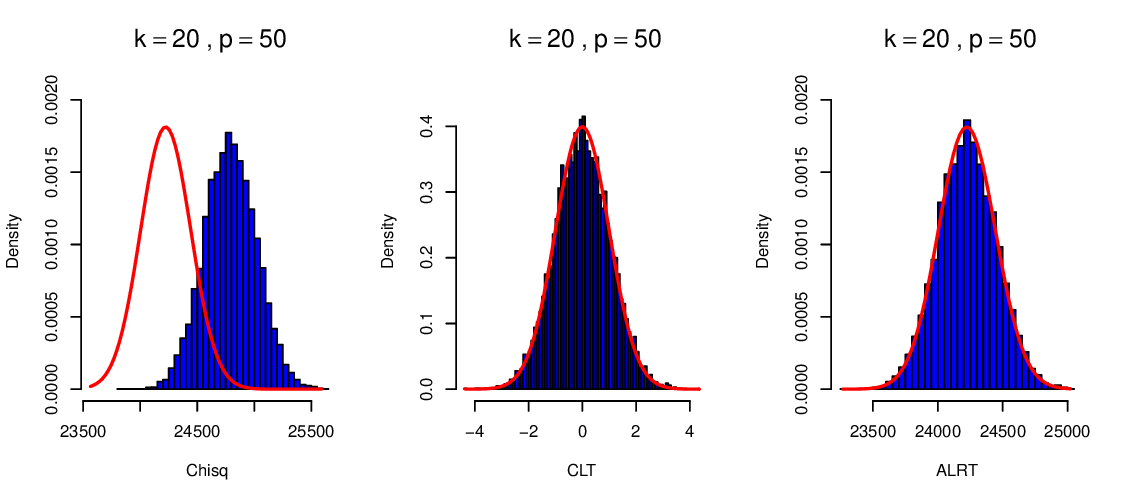}
\includegraphics[height=2.0in,width=6in]{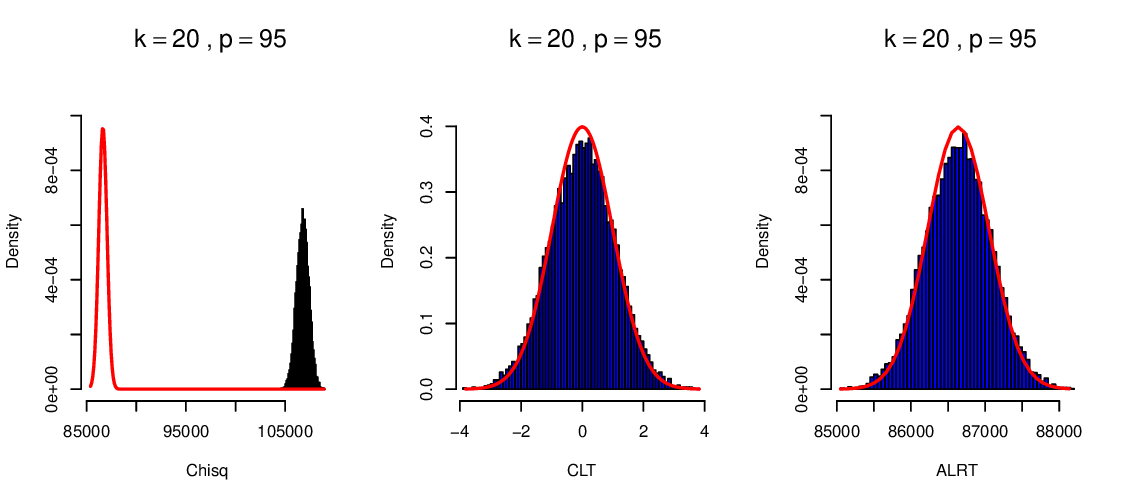}
\caption{Histograms of $-2\rho \log\Lambda_n^*$ (Chisq),
$(-2\log\Lambda_n^*-\mu_n)/\sigma_n$ (CLT), and  $Z_n$ (ALRT), where
$n_i=100$, $1 \le i\le k$, with $k=20$ and $p=5,20,50,95$.}
 \label{histogram-k20}
\end{figure}

\begin{figure}[H]
\centering
\includegraphics[height=2.0in,width=6in]{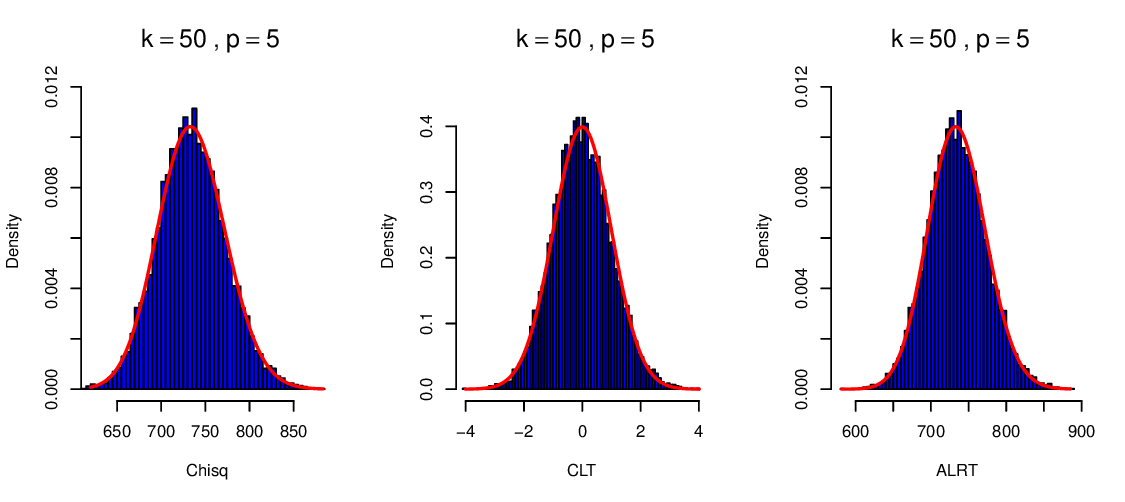}
\includegraphics[height=2.0in,width=6in]{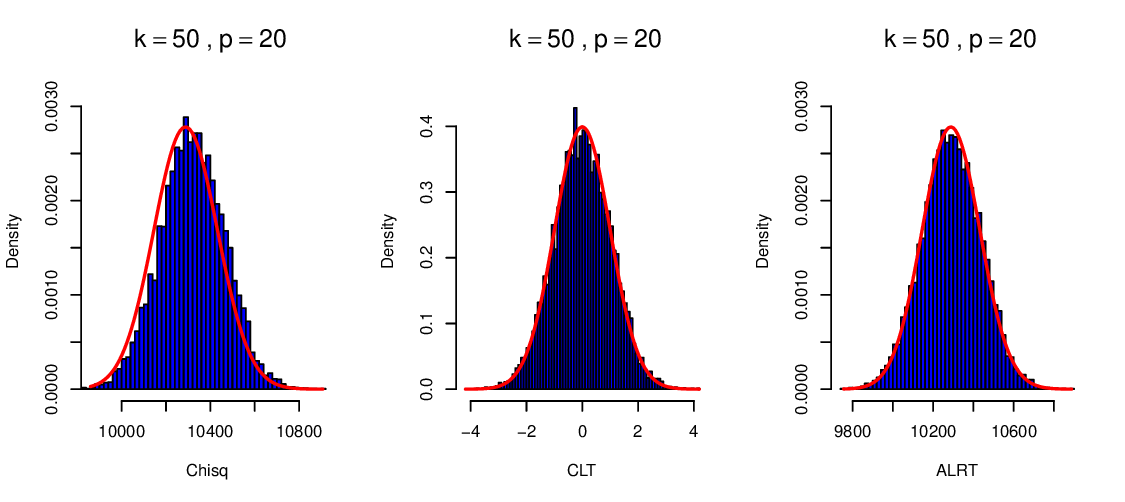}
\includegraphics[height=2.0in,width=6in]{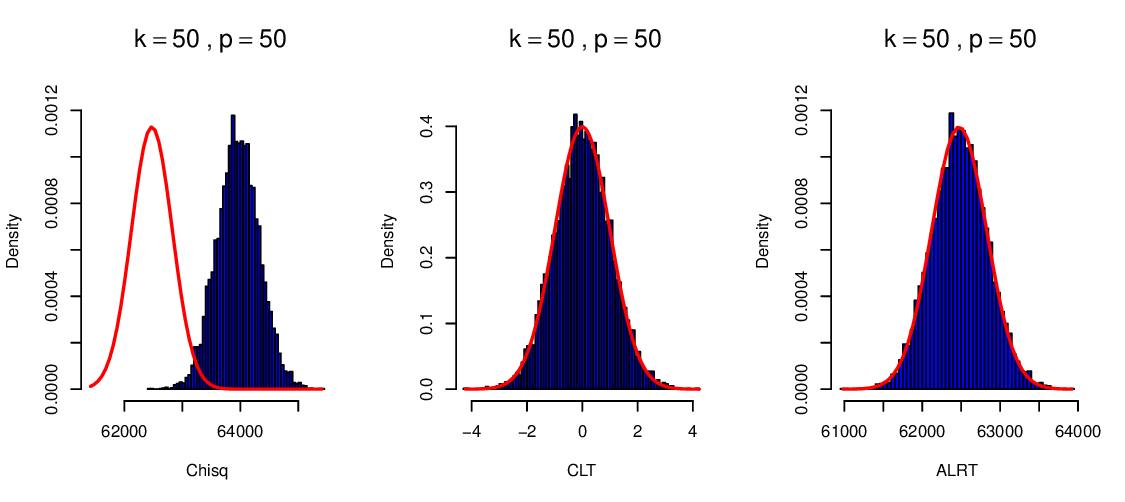}
\includegraphics[height=2.0in,width=6in]{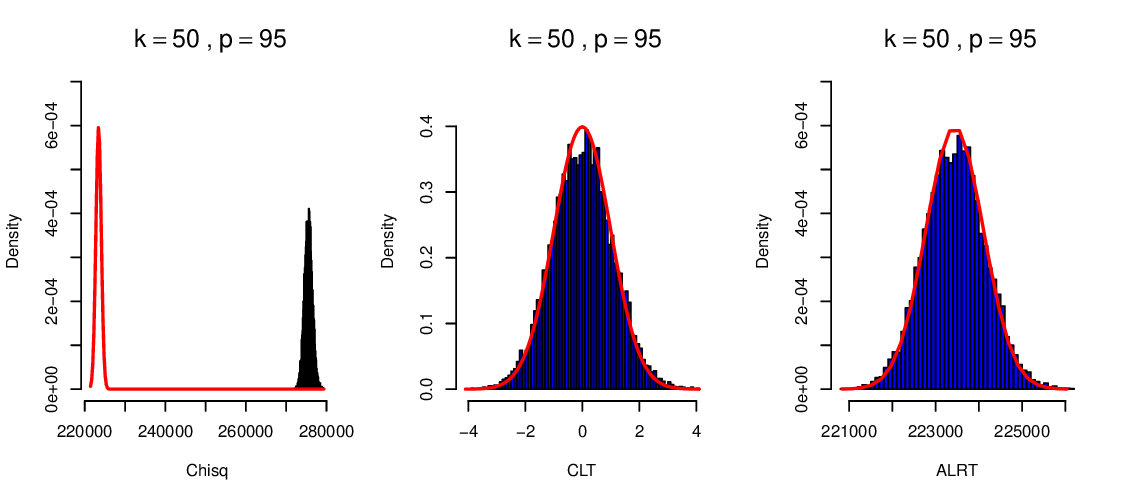}
\caption{Histograms of $-2\rho \log\Lambda_n^*$ (Chisq),
$(-2\log\Lambda_n^*-\mu_n)/\sigma_n$ (CLT), and  $Z_n$ (ALRT), where
$n_i=100$, $1 \le i\le k$, with $k=50$ and $p=5,20,50,95$.}
\label{histogram-k50}
\end{figure}

\end{document}